\documentclass[12pt]{article}
\usepackage[a4paper,margin=1in]{geometry}
\usepackage{amsmath,amssymb,amsthm,mathtools,bm}
\usepackage{graphicx,booktabs,array,enumitem}
\usepackage{xcolor,hyperref,microtype,authblk,caption,subcaption,natbib,placeins}
\usepackage{xspace}
\usepackage{float}
\usepackage{etoolbox}
\floatstyle{ruled}
\newfloat{algorithm}{tbp}{loa}
\floatname{algorithm}{Algorithm}
\usepackage{cleveref}
\crefname{algorithm}{Algorithm}{Algorithms}
\Crefname{algorithm}{Algorithm}{Algorithms}
\graphicspath{{figures/}}
\hypersetup{colorlinks=true,linkcolor=blue!55!black,citecolor=blue!55!black,urlcolor=blue!55!black}
\numberwithin{equation}{section}
\newtheorem{theorem}{Theorem}[section]
\newtheorem{lemma}[theorem]{Lemma}

\newtheorem{corollary}[theorem]{Corollary}
\newtheorem{assumption}[theorem]{Assumption}
\crefname{assumption}{Assumption}{Assumptions}
\Crefname{assumption}{Assumption}{Assumptions}
\theoremstyle{remark}
\newtheorem{remark}[theorem]{Remark}

\crefname{lemma}{Lemma}{Lemmas}
\Crefname{lemma}{Lemma}{Lemmas}
\AtBeginEnvironment{lemma}{\crefalias{theorem}{lemma}}
\AtEndEnvironment{lemma}{\crefalias{theorem}{theorem}}
\AtBeginEnvironment{assumption}{\crefalias{theorem}{assumption}}
\AtEndEnvironment{assumption}{\crefalias{theorem}{theorem}}
\newcommand{\X}{\mathcal X}
\newcommand{\Y}{\mathcal Y}
\newcommand{\B}{\mathcal B}
\newcommand{\eps}{\varepsilon}
\newcommand{\norm}[1]{\left\lVert #1\right\rVert}
\newcommand{\ip}[2]{\left\langle #1,#2\right\rangle}
\newcommand{\R}{\mathbb R}

\newcommand{\RG}{RG-AA-LW\xspace}

\title{\textbf{Anderson-Accelerated Landweber Iteration with Regularization Guarantees for Nonlinear Ill-Posed Problems}}
\makeatletter
\renewcommand{\@author}{\normalfont\fontsize{12}{14}\selectfont
Jianming Gao\textsuperscript{1}, Xiliang Lu\textsuperscript{1,2,3}, and Lin Tian\textsuperscript{1,*}}
\patchcmd{\@maketitle}{\vskip 1em}{\vskip 0pt}{}{\PackageError{titlelayout}{Date-space patch failed}{}}
\patchcmd{\@maketitle}{\par\vskip 1.5em}{\par\vskip 0pt}{}{\PackageError{titlelayout}{Bottom-space patch failed}{}}
\makeatother
\date{}

\begin{document}
\maketitle
\begingroup
\makeatletter
\renewcommand{\@makefntext}[1]{\noindent#1}
\makeatother
\renewcommand{\thefootnote}{}
\hypersetup{urlcolor=black}
\footnotetext{\normalfont\fontsize{10}{12}\selectfont%
\setlength{\parskip}{2pt}%
\noindent\hspace*{1em}\begin{tabular}{@{}l@{}}
Jianming Gao\\
\href{mailto:2020302011052@whu.edu.cn}{2020302011052@whu.edu.cn}\\[3pt]
Xiliang Lu\\
\href{mailto:xllv.math@whu.edu.cn}{xllv.math@whu.edu.cn}\\[3pt]
Lin Tian\\
\href{mailto:tianlin0000000@163.com}{tianlin0000000@163.com}
\end{tabular}\par
\vspace{5pt}
\hangindent=1em\hangafter=1\noindent\makebox[1em][l]{\textsuperscript{1}}School of Mathematics and Statistics, Wuhan University, Wuhan 430072, China\par
\hangindent=1em\hangafter=1\noindent\makebox[1em][l]{\textsuperscript{2}}Hubei Center for Applied Mathematics, Wuhan University, Wuhan 430072, China\par
\hangindent=1em\hangafter=1\noindent\makebox[1em][l]{\textsuperscript{3}}Hubei Computational Science Key Laboratory, Wuhan University, Wuhan 430072, China\par
\noindent\textsuperscript{*}Corresponding author.
}
\endgroup
\vspace*{-30pt}

\begin{abstract}
Anderson acceleration (AA) can markedly reduce iteration counts for fixed-point solvers, but its unrestricted signed history extrapolation need not preserve the stability required for iterative regularization of ill-posed problems.  We introduce a regularization-guaranteed Anderson-accelerated Landweber method (\RG) in which finite-memory AA generates a proposal while the committed iterate is a controlled perturbation of a certified Landweber anchor.  Under a local tangential cone condition and a scale-normalized localization condition, we prove strong convergence of the exact-data trajectory to a solution, stability at every fixed iteration, finite discrepancy stopping, and strong convergence of the discrepancy-stopped noisy reconstructions to the same exact-data limit; local uniqueness is needed only if that selected limit is to be identified with a prescribed solution $x^\dagger$.  Under the additional benchmark source condition
$x^\dagger-x_0=F'(x^\dagger)^*\lambda^\dagger$ and a standard range-invariance condition, we obtain the discrepancy-principle rate
$\norm{x_{k_\delta}^\delta-x^\dagger}=O(\delta^{1/2})$.  The rate proof controls data-dependent Anderson coefficients through fixed memory, a fixed-reference ridge in the small AA least-squares problem, and an outer correction budget with an $O(k^{-1})$ tail. The numerical examples verify the efficiency of the proposed method. 
\end{abstract}


\section{Introduction}
Let $F:D(F)\subset\X\to\Y$ be a nonlinear operator between real Hilbert spaces and consider
\begin{equation}\label{eq:inverse}
 F(x)=y,
\end{equation}
where only noisy data $y^\delta$ satisfying
\begin{equation}\label{eq:noise}
 \norm{y^\delta-y}\le\delta
\end{equation}
are available.  Because the inverse map is generally discontinuous, small perturbations of the data may produce large changes in a solution.  A convergent iterative solver is therefore not sufficient: the iteration index must act as a regularization parameter and must be chosen in a noise-dependent manner \citep{EnglHankeNeubauer1996,KNS2008}.  Nonlinear Landweber iteration realizes this principle using only evaluations of $F$ and $F'(x)^*$.  Under a tangential cone condition, its local convergence, stability at fixed indices, finite Morozov stopping, and iterative-regularization property are well understood \citep{HNS1995,KNS2008,Neubauer2000}.

Several acceleration strategies preserve this regularization structure.  Linear semi-iterative methods achieve order-optimal rates under suitable source conditions \citep{Hanke1991}.  For nonlinear problems, Nesterov-type, two-point-gradient, and adaptive heavy-ball methods combine extrapolation with noise-dependent stopping and convergence analysis \citep{HubmerRamlau2018,ZhongWangJin2019,JinHuang2024}.  Regularizing Newton and Levenberg--Marquardt methods can reduce the number of outer iterations further, but require the approximate solution of a regularized linearized inverse problem at each step \citep{Hanke1997,KNS2008}.

Anderson acceleration (AA) uses a different mechanism.  It forms a signed affine combination of recent fixed-point residuals by solving a small least-squares problem \citep{Anderson1965,FangSaad2009}.  For linear problems, full-memory AA is closely related to GMRES; for nonlinear problems, it can be interpreted as a multisecant update \citep{WalkerNi2011,PotraEngler2013}.  These viewpoints explain why AA can remove slowly decaying components much faster than the underlying fixed-point iteration.  Local convergence analyses confirm this effect for contractive and certain noncontractive maps \citep{TothKelley2015,EvansEtAl2020,PollockRebholz2021}.  They also identify a basic instability: nearly dependent histories can produce large signed coefficients and unreliable extrapolation.  Regularized coefficient problems, restarting, and globalization have consequently been developed for well-posed fixed-point computations \citep{ScieurEtAl2016,GarstkaEtAl2022,OuyangEtAl2023}.

In an inverse problem, rapid residual reduction and regularization are distinct requirements.  Even in the linear setting, GMRES-type methods need early stopping or an additional regularization device when applied to noisy ill-posed systems \citep{Reichel2002}.  The AA least-squares residual alone controls neither the coefficient vector nor the distance of the extrapolated point from a solution.  Consequently, unrestricted AA applied to Landweber iteration may accelerate both recoverable solution components and noise-dominated components.  A residual-only acceptance test does not recover the stability estimates on which iterative regularization relies.

AA has shown considerable computational promise in inverse problems.  \citet{JinZhou2021} accelerated a fixed-point reconstruction for a subdiffusion potential.  \citet{Yang2021} applied AA to steepest-descent seismic inversion, and \citet{AghazadeEtAl2022} accelerated an augmented-Lagrangian formulation of extended waveform inversion.  AA has also been incorporated into learned EIT reconstruction \citep{XuWangZhou2024} and an adaptively damped variational Born iteration for electromagnetic inverse scattering \citep{ZhouEtAl2025}.  Most closely related to the present work, \citet{ZhangHeFu2026} developed an AA-based fast Landweber-type method for nonlinear ill-posed problems and demonstrated substantial numerical gains.  The present paper addresses the complementary question that is decisive for ill-posed problems: how much of the AA proposal can be committed without destroying the stability mechanism of iterative regularization.  We retain a standard finite-memory AA proposal but place the regularization guarantee in a separate commitment layer.

The method separates proposal from commitment.  A standard finite-memory AA step first generates a candidate.  The accepted iterate is then obtained by radially truncating the difference between this candidate and a certified Landweber anchor.  The admissible radius is the smaller of a summable global budget and a fixed multiple of the current Landweber displacement.  The global budget makes the accepted trajectory a summably perturbed Landweber sequence, while the relative gate forces the perturbation to vanish with the certified step.  A fixed-reference ridge uniquely defines the small AA coefficient problem and yields finite-index continuity.  The Landweber anchor, safeguard, and discrepancy principle regularize the inverse problem.  When the raw AA proposal lies inside the admissible ball, it is accepted unchanged.

This construction leads to two levels of theory.  Under a local tangential cone condition and a scale-normalized localization condition, the exact-data iterates converge strongly to a solution $x_*$, the noisy iterates are stable at every fixed index, the discrepancy principle stops after finitely many steps, and the stopped reconstructions converge strongly to the same selected limit $x_*$ as $\delta\to0$.  Thus local uniqueness is not required for the regularization statement itself; it is only one way to identify $x_*$ with a prescribed solution $x^\dagger$.  A convergence rate toward $x^\dagger$ requires additional structure.  We impose the range-invariance condition used by \citet{Jin2025},
\[
 F'(x)=Q_xF'(x^\dagger),\qquad
 \|I-Q_x\|\le \kappa_0\|x-x^\dagger\|,
\]
together with the benchmark source condition $x^\dagger-x_0=F'(x^\dagger)^*\lambda^\dagger$.  Following the source-element strategy of \citet{Jin2025}, we show that the AA history produces an additional source-space perturbation.  Fixed memory and the fixed-reference ridge bound this perturbation by a finite-window residual energy, so no a priori uniform bound on the Anderson coefficients is required.  If the summable outer budget has an $O(k^{-1})$ tail, the discrepancy-stopped error is $O(\sqrt\delta)$.

The main contributions of this paper are as follows:
\begin{enumerate}[label=(\roman*),leftmargin=2.2em]
\item We formulate \RG as a standard finite-memory AA proposal followed by a two-scale radial commitment rule around a certified Landweber anchor.  The algorithm uses only a small regularized least-squares problem in addition to the evaluations required by the underlying first-order iteration.
\item We prove strong exact-data convergence, fixed-index stability, finite discrepancy stopping, and the iterative-regularization property without assuming local uniqueness: the noisy discrepancy-stopped reconstructions converge to the same solution selected by the exact-data trajectory.  The proof reduces the accepted AA correction to a summable perturbation whose size also vanishes with the Landweber displacement.
\item Under the additional benchmark source and range-invariance conditions, we prove source-element stability through $O(\delta^{-1})$ iterations and obtain
$\norm{x_{k_\delta}^\delta-x^\dagger}=O(\delta^{1/2})$ whenever
$\sup_{n\ge1}n\sum_{k=n}^\infty\omega_k<\infty$.
\end{enumerate}

The numerical section explains three questions: verification of the theorem-aligned rate, robustness outside the verified structural assumptions, and computational benefit when the AA proposal is compatible with the certified trajectory.  In difficult regimes the safeguard may trade some raw-AA speed for stability; when the proposal is already admissible, it is accepted unchanged.

The paper is organized as follows.  Section~2 introduces Landweber iteration, raw AA--Landweber, and the safeguarded \RG method in that order.  Section~3 establishes the regularization properties.  Section~4 proves the benchmark convergence rate.  Sections~5 and~6 present the numerical experiments and conclusions.

\section{Proposed \RG algorithm}

Throughout, $\X$ and $\Y$ are real Hilbert spaces.  The method is stopped at the first index $k$ for which
\begin{equation}\label{eq:discrepancy}
 \norm{F(x_k^\delta)-y^\delta}\le\tau\delta,
\end{equation}
where $\tau>1$ is fixed.  We describe the three ingredients of the method in increasing order of complexity.

\subsection{Landweber iteration and the certified anchor}

The nonlinear Landweber iteration with a fixed step length $\mu>0$ is
\begin{equation}\label{eq:landweber}
 x_{k+1}^\delta
 =x_k^\delta-\mu F'(x_k^\delta)^*\bigl(F(x_k^\delta)-y^\delta\bigr).
\end{equation}
For the safeguarded method, the step length is chosen from the current residual.  Set
\begin{equation}\label{eq:rga}
 r_k^\delta=F(x_k^\delta)-y^\delta,
 \qquad g_k^\delta=F'(x_k^\delta)^*r_k^\delta,
\end{equation}
and fix $0\le\eta<1$, $0<\bar\gamma<2$, and
\begin{equation}\label{eq:tau}
 \tau>\frac{1+\eta}{1-\eta},
 \qquad c_\tau:=1-\eta-\frac{1+\eta}{\tau}>0.
\end{equation}
Before \eqref{eq:discrepancy} is satisfied, define
\begin{equation}\label{eq:adef}
 a_k^\delta=(1-\eta)\norm{r_k^\delta}^2
 -(1+\eta)\delta\norm{r_k^\delta},
\end{equation}
\begin{equation}\label{eq:tdef}
 t_k^\delta
 =\min\left\{\bar t,
 \bar\gamma\frac{a_k^\delta}{\norm{g_k^\delta}^2}\right\},
 \qquad \bar t>0,
 \qquad
 z_{k+1}^\delta=x_k^\delta-t_k^\delta g_k^\delta.
\end{equation}
The point $z_{k+1}^\delta$ is the certified Landweber anchor.  The analysis in Section~3 shows that $a_k^\delta>0$ and $g_k^\delta\ne0$ before stopping, so \eqref{eq:tdef} is well defined.  Notice that $t_k^\delta$ is chosen solely to certify descent of the committed trajectory.  It is distinct from the fixed step $\mu$ introduced below, which belongs only to the raw AA proposal map.

\subsection{Raw finite-memory AA--Landweber}
For the same fixed step $\mu>0$, define
\begin{equation}\label{eq:rawmap}
 G_\mu^\delta(x)=x-\mu F'(x)^*(F(x)-y^\delta).
\end{equation}

For committed iterates set
\[
 f_j^\delta:=G_\mu^\delta(x_j^\delta)-x_j^\delta.
\]
Let $m\ge1$ be fixed.  At iteration $k$, let $\ell_k=\min\{m,k\}$ and use the history $j=k-\ell_k,\ldots,k$.  Form the differences
\[
 \Delta F_k^\delta=[f_{j+1}^\delta-f_j^\delta]_{j=k-\ell_k}^{k-1},
 \qquad
 \Delta G_k^\delta=[G_\mu^\delta(x_{j+1}^\delta)-G_\mu^\delta(x_j^\delta)]_{j=k-\ell_k}^{k-1}.
\]
When $k=0$ these matrices have zero columns.  Once $f_0^\delta$ is available, fix the reference-scale ridge
\begin{equation}\label{eq:ridge}
 \lambda_{\rm AA}^\delta:=\eps_{\rm AA}\norm{f_0^\delta}^2,
 \qquad \eps_{\rm AA}>0.
\end{equation}
Setting $\eps_{\rm AA}=0$ gives the usual raw AA coefficient problem.  The proposed method uses $\eps_{\rm AA}>0$, which makes the proposal unique and continuous with respect to every finite history.
For $k\ge1$ define
\begin{equation}\label{eq:aa-LS}
 \gamma_k^\delta=\arg\min_{\gamma\in\R^{\ell_k}}
 \left\{
 \norm{f_k^\delta-\Delta F_k^\delta\gamma}^2
 +\lambda_{\rm AA}^\delta\norm{\gamma}_2^2
 \right\}.
\end{equation}
The raw AA proposal is
\begin{equation}\label{eq:aa-cand}
 \widehat{x}_{k+1}^{\rm AA,\delta}
 =G_\mu^\delta(x_k^\delta)-\Delta G_k^\delta\gamma_k^\delta,
\end{equation}
with $\widehat{x}_1^{\rm AA,\delta}=G_\mu^\delta(x_0)$.
Committing this proposal directly, $x_{k+1}^\delta=\widehat{x}_{k+1}^{\rm AA,\delta}$, gives raw AA--Landweber.  The proposed method changes only this commitment step.

\begin{remark}
The difference form \eqref{eq:aa-LS}--\eqref{eq:aa-cand} is algebraically equivalent to the standard affine AA formulation when the least-squares problem is unregularized.  We use this form because it exposes the finite-memory history correction in the source-element analysis.  The ridge stabilizes and uniquely defines the finite-dimensional Anderson coefficient problem.  The Landweber anchor, outer safeguard, and discrepancy stopping rule provide inverse-problem regularization.
\end{remark}

\begin{remark}[Two Landweber roles]
The fixed step $\mu$ and the certified step $t_k^\delta$ need not coincide.  The map $G_\mu^\delta$ is used to build the AA history and hence determines the proposed acceleration direction.  By contrast, $t_k^\delta$ is chosen from the current residual so that the anchor satisfies the descent estimate in \Cref{lem:basicinner}.  The algorithm never certifies the raw map itself; it certifies only the committed perturbation of the anchor.
\end{remark}

\subsection{The two-scale \RG safeguard}
Let $\{\omega_k\}_{k\ge0}$ be deterministic and nonnegative with
\begin{equation}\label{eq:omegasum}
 \sum_{k=0}^\infty\omega_k=\Omega<\infty.
\end{equation}
Define the certified displacement
\begin{equation}\label{eq:sigma}
 \sigma_k^\delta:=\norm{z_{k+1}^\delta-x_k^\delta}
 =t_k^\delta\norm{g_k^\delta}.
\end{equation}
Fix $\chi>0$ and set
\begin{equation}\label{eq:rho}
 \rho_k^\delta=\min\{\omega_k\sigma_0^\delta,\;\chi\sigma_k^\delta\}.
\end{equation}
\begin{equation}\label{eq:correction}
 x_{k+1}^\delta
 =z_{k+1}^\delta
 +P_{\rho_k^\delta\B}\bigl(\widehat{x}_{k+1}^{\rm AA,\delta}-z_{k+1}^\delta\bigr).
\end{equation}
where $\B:=\{v\in\X:\norm v\le1\}$ and $P_{\rho\B}$ denotes radial projection onto the ball of radius $\rho$ centered at the origin.  Equivalently,
\[
 e_k^\delta:=x_{k+1}^\delta-z_{k+1}^\delta
 =\beta_k^\delta\bigl(\widehat{x}_{k+1}^{\rm AA,\delta}-z_{k+1}^\delta\bigr),
 \hspace{0.8em}
 \beta_k^\delta=
 \begin{cases}
 \min\{1,\rho_k^\delta/\norm{\widehat{x}_{k+1}^{\rm AA,\delta}-z_{k+1}^\delta}\},
 &\widehat{x}_{k+1}^{\rm AA,\delta}\ne z_{k+1}^\delta,\\
 0,&\widehat{x}_{k+1}^{\rm AA,\delta}=z_{k+1}^\delta.
 \end{cases}
\]

\begin{algorithm}[t]
\caption{Regularization-guaranteed AA--Landweber iteration (\RG)}
\label{alg:rglAA}
\begin{enumerate}[leftmargin=2.5em,label=\arabic*.]
\item \textbf{Input:} initial guess $x_0$, noise level $\delta$, discrepancy
parameter $\tau$, memory depth $m$, raw Landweber step $\mu$, certified-step
parameters $(\eta,\bar\gamma,\bar t)$, ridge parameter $\eps_{\rm AA}>0$,
correction schedule $\{\omega_k\}_{k\ge0}$, and relative safeguard parameter
$\chi>0$.  Set $x_0^\delta=x_0$.
\item For $k=0,1,2,\ldots$, perform the following steps:
  \begin{enumerate}[leftmargin=2.4em,label=(\alph*)]
  \item Set $r_k^\delta=F(x_k^\delta)-y^\delta$.  If
  $\norm{r_k^\delta}\le\tau\delta$, set $k_\delta=k$ and return
  $x_{k_\delta}^\delta$.
  \item Compute $g_k^\delta=F'(x_k^\delta)^*r_k^\delta$, determine
  $t_k^\delta$ from \eqref{eq:tdef}, and set
  $z_{k+1}^\delta=x_k^\delta-t_k^\delta g_k^\delta$ and
  $\sigma_k^\delta=\norm{z_{k+1}^\delta-x_k^\delta}$.
  \item Compute
  $f_k^\delta=G_\mu^\delta(x_k^\delta)-x_k^\delta$.
  If $k=0$, set
  $\lambda_{\rm AA}^\delta=\eps_{\rm AA}\norm{f_0^\delta}^2$ and store
  $\sigma_0^\delta$.
  \item Compute $\widehat{x}_{k+1}^{\rm AA,\delta}$ from
  \eqref{eq:aa-LS}--\eqref{eq:aa-cand}.
  \item Set
  $\rho_k^\delta=\min\{\omega_k\sigma_0^\delta,\chi\sigma_k^\delta\}$.
  \item Update
  \[
    x_{k+1}^\delta=z_{k+1}^\delta+
    P_{\rho_k^\delta\B}
    \bigl(\widehat{x}_{k+1}^{\rm AA,\delta}-z_{k+1}^\delta\bigr).
  \]
  \end{enumerate}
\end{enumerate}
\end{algorithm}

\begin{remark}[Role of the tangential-cone bound in implementation]
The theory assumes that a valid upper bound $\eta$ for the local tangential-cone constant is available, because $\eta$ enters both the certified anchor and the admissible discrepancy parameter.  In applications this constant is rarely known sharply.  A conservative bound may be used when problem structure provides one; in the experiments below we use fixed values of $\eta$.  Designing a fully observable or backtracking-type certification that removes prior knowledge of $\eta$ is left for future work.
\end{remark}

If an exact or noisy trajectory satisfies its stopping test, we extend it
constantly for all subsequent indices.  This convention is used only when
comparing fixed indices across noise levels; the reported stopping index and
the algorithmic output are unchanged.  The case $F(x_0)=y$ is trivial for
exact data and terminates at index zero.  All coefficient estimates below are
therefore stated for the nontrivial case $F(x_0)\ne y$.

\begin{remark}[Rate-compatible schedules]
The regularization results require only the summability condition \eqref{eq:omegasum}.  A rate-compatible example is
\begin{equation}\label{eq:poweromega}
 \omega_k=\frac{\Omega}{\zeta(s)}(k+1)^{-s},
 \qquad s\ge2,
\end{equation}
where
\[
 \zeta(s):=\sum_{j=1}^{\infty}j^{-s}
\]
is the Riemann zeta function.  This normalization gives
$\sum_{k=0}^{\infty}\omega_k=\Omega$.  Section~4 states the additional tail
condition together with the other rate assumptions.
\end{remark}

\section{Regularization properties}

This section collects all assumptions used for convergence and regularization before any auxiliary result is stated.  We write $x_k=x_k^0$, $r_k=r_k^0$, and so forth for exact data.  The case $F(x_0)=y$ terminates at $k=0$; hence only the nontrivial case $F(x_0)\ne y$ is considered below.

\subsection{Assumptions}

\begin{assumption}[Local tangential cone condition]\label{ass:tcc}
Let $x^\dagger$ be a solution of $F(x)=y$.  There exist $R>0$, $0\le\eta<1$, and $M>0$ such that $\overline B_R(x^\dagger)\subset D(F)$, the operator $F$ is continuously Fr\'echet differentiable on this ball, and
\begin{equation}\label{eq:tcc}
 \norm{F(z)-F(x)-F'(x)(z-x)}
 \le \eta\norm{F(z)-F(x)},
 \qquad x,z\in\overline B_R(x^\dagger),
\end{equation}
\begin{equation}\label{eq:derivativebound}
 \norm{F'(x)}\le M,
 \qquad x\in\overline B_R(x^\dagger).
\end{equation}
The constant $\eta$ is the one used in \eqref{eq:adef}, and $\tau$ satisfies \eqref{eq:tau}.
\end{assumption}

For exact data, let
\[
 d_0:=\norm{x_0-x^\dagger},
 \qquad \sigma_0:=\norm{z_1-x_0}.
\]

\begin{assumption}[Scale-normalized localization]\label{ass:margin}
The summable budget in \eqref{eq:omegasum} satisfies
\begin{equation}\label{eq:margin}
 d_0+\Omega\sigma_0<R.
\end{equation}
Fix numbers $\bar\sigma_0>\sigma_0$ and $R_*<R$ such that
\begin{equation}\label{eq:Rstar}
 d_0+\Omega\bar\sigma_0\le R_*<R.
\end{equation}
The number $\bar\sigma_0$ provides a uniform upper bound for the first noisy
certified displacement $\sigma_0^\delta$ when $\delta$ is sufficiently small,
and $R_*$ is the common localization radius for all corresponding noisy
trajectories.
\end{assumption}

Assumptions~\ref{ass:tcc} and~\ref{ass:margin} are the only assumptions used in Section~3.  They imply convergence to the solution selected by the exact-data trajectory; no local uniqueness hypothesis is needed for this regularization statement.  If the local solution is unique, the selected limit is of course the prescribed solution $x^\dagger$.

\begin{remark}[Interpretation of the localization condition]
Condition \eqref{eq:margin} is a sufficient local-theory assumption rather than an a priori computable stopping rule: it is written relative to the unknown solution and guarantees that the summably perturbed trajectory remains inside the neighborhood where the TCC is valid.  The scale $\sigma_0$ makes the global AA budget dimensionally tied to the certified dynamics.  A posteriori localization tests or adaptive correction budgets would be preferable in applications and are discussed again in Section~6.
\end{remark}

\subsection{Preliminary estimates}

Before discrepancy stopping, \eqref{eq:tau} gives
\begin{equation}\label{eq:a-lower}
 a_k^\delta\ge c_\tau\norm{r_k^\delta}^2>0.
\end{equation}

\begin{lemma}[Certified Landweber descent]\label{lem:basicinner}
Under \Cref{ass:tcc}, every pre-discrepancy iterate in $\overline B_R(x^\dagger)$ satisfies
\begin{equation}\label{eq:innerlower}
 \ip{g_k^\delta}{x_k^\delta-x^\dagger}\ge a_k^\delta>0.
\end{equation}
Consequently, \eqref{eq:tdef} is well defined and
\begin{equation}\label{eq:t-lower}
 t_k^\delta\ge t_*:=\min\left\{\bar t,\frac{\bar\gamma c_\tau}{M^2}\right\}>0,
\end{equation}
\begin{equation}\label{eq:safe-rdescent}
 \norm{z_{k+1}^\delta-x^\dagger}^2
 \le \norm{x_k^\delta-x^\dagger}^2
 -\kappa_\tau\norm{r_k^\delta}^2,
 \qquad \kappa_\tau:=(2-\bar\gamma)t_*c_\tau>0.
\end{equation}
\end{lemma}

\begin{proof}
Put $r=F(x_k^\delta)-y$.  Applying \eqref{eq:tcc} with $x=x_k^\delta$ and $z=x^\dagger$ yields
\[
 \norm{F'(x_k^\delta)(x_k^\delta-x^\dagger)-r}\le\eta\norm r.
\]
Since $r=r_k^\delta+(y^\delta-y)$ and $\norm r\le\norm{r_k^\delta}+\delta$, we obtain
\begin{align*}
 \ip{g_k^\delta}{x_k^\delta-x^\dagger}
 &\ge \norm{r_k^\delta}^2-\delta\norm{r_k^\delta}
       -\eta\norm r\,\norm{r_k^\delta}\\
 &\ge (1-\eta)\norm{r_k^\delta}^2
       -(1+\eta)\delta\norm{r_k^\delta}=a_k^\delta.
\end{align*}
This proves \eqref{eq:innerlower}, and therefore $g_k^\delta\ne0$.  Expanding the square in \eqref{eq:tdef} and using
$t_k^\delta\norm{g_k^\delta}^2\le\bar\gamma a_k^\delta$ gives
\[
 \norm{z_{k+1}^\delta-x^\dagger}^2
 \le\norm{x_k^\delta-x^\dagger}^2
 -(2-\bar\gamma)t_k^\delta a_k^\delta.
\]
Finally, $\norm{g_k^\delta}\le M\norm{r_k^\delta}$ and \eqref{eq:a-lower} imply \eqref{eq:t-lower}; combining these bounds proves \eqref{eq:safe-rdescent}.
\end{proof}

Because the first certified displacement depends continuously on the data, there is a $\delta_0>0$ such that $\sigma_0^\delta\le\bar\sigma_0$ whenever $0\le\delta\le\delta_0$.  For later use, set
\begin{equation}\label{eq:bareps}
 \bar\eps_k:=\bar\sigma_0\omega_k,
 \qquad \xi_k:=2R_*\bar\eps_k+\bar\eps_k^2,
 \qquad \theta:=\sum_{k=0}^\infty\xi_k<\infty.
\end{equation}

\begin{lemma}[Localization and quasi-Fej\'er estimate]\label{lem:qf}
Let \Cref{ass:tcc,ass:margin} hold.  For $0\le\delta\le\delta_0$, every committed iterate before discrepancy stopping satisfies
\begin{equation}\label{eq:local}
 \norm{x_k^\delta-x^\dagger}\le d_0+\Omega\sigma_0^\delta\le R_*<R,
\end{equation}
and
\begin{equation}\label{eq:qf}
 \norm{x_{k+1}^\delta-x^\dagger}^2
 \le \norm{x_k^\delta-x^\dagger}^2
 -\kappa_\tau\norm{r_k^\delta}^2+\xi_k.
\end{equation}
Consequently, for every $N\le k_\delta$,
\begin{equation}\label{eq:resenergy}
 \sum_{k=0}^{N-1}\norm{r_k^\delta}^2\le C_r,
\end{equation}
where $C_r$ is independent of $N$ and sufficiently small $\delta$.
\end{lemma}
\begin{proof}
The safe anchor is Fej\'er decreasing by \Cref{lem:basicinner}.  Since $\norm{e_k^\delta}\le\omega_k\sigma_0^\delta$, induction gives
\[
 \norm{x_k^\delta-x^\dagger}
 \le d_0+\sigma_0^\delta\sum_{j=0}^{k-1}\omega_j
 \le R_*.
\]
Thus every step remains in the ball on which \Cref{lem:basicinner} applies.  Expanding $x_{k+1}^\delta=z_{k+1}^\delta+e_k^\delta$, using $\norm{z_{k+1}^\delta-x^\dagger}\le R_*$, and then $\norm{e_k^\delta}\le\bar\eps_k$ proves \eqref{eq:qf}.  Summation gives \eqref{eq:resenergy} with $C_r=(d_0^2+\theta)/\kappa_\tau$.
\end{proof}

For exact data, the same estimate gives
\begin{equation}\label{eq:exactqf}
 \norm{x_{k+1}-x^\dagger}^2
 \le \norm{x_k-x^\dagger}^2-c_0t_k\norm{r_k}^2+\xi_k,
 \qquad c_0=(2-\bar\gamma)(1-\eta)>0,
\end{equation}
where $\xi_k$ is defined in \eqref{eq:bareps}.  Since $t_k$ has a uniform
positive lower bound and $\sum_k\xi_k<\infty$, it follows that
$\sum_k\norm{r_k}^2<\infty$ and $\norm{x_k-x^\dagger}^2$ has a finite limit.

\begin{lemma}[TCC cross estimate]\label{lem:cross}
Let $j>k$.  Choose $\ell\in\{k,\ldots,j\}$ so that
$\norm{r_\ell}=\min_{k\le i\le j}\norm{r_i}$.  Then, for every $i\in\{k,\ldots,j\}$,
\begin{equation}\label{eq:cross}
 \norm{F'(x_i)(x_\ell-x^\dagger)}\le C_\eta\norm{r_i},
\end{equation}
where $C_\eta$ depends only on the tangential cone constant.
\end{lemma}
\begin{proof}
Apply \eqref{eq:tcc} first to $(x_i,x_\ell)$ and then to $(x_i,x^\dagger)$.  This gives
\begin{align*}
 \norm{F'(x_i)(x_\ell-x_i)}
 &\le(1+\eta)\bigl(\norm{r_\ell}+\norm{r_i}\bigr),\\
 \norm{F'(x_i)(x_i-x^\dagger)}
 &\le(1+\eta)\norm{r_i}.
\end{align*}
Since $\norm{r_\ell}\le\norm{r_i}$, the triangle inequality proves \eqref{eq:cross} with $C_\eta=3(1+\eta)$.
\end{proof}

\subsection{Strong convergence for exact data}

\begin{theorem}[Strong exact-data convergence]\label{thm:strong}
Under \Cref{ass:tcc,ass:margin}, the exact-data iterates generated by Algorithm~\ref{alg:rglAA} form a Cauchy sequence in $\X$.  Hence $x_k\to x_*$ strongly and $F(x_*)=y$.  If $x^\dagger$ is the unique solution in $\overline B_{R_*}(x^\dagger)$, then $x_*=x^\dagger$.
\end{theorem}
\begin{proof}
The committed increment has the exact form
\[
 x_{i+1}-x_i=-t_iF'(x_i)^*r_i+e_i,
 \qquad \norm{e_i}\le\bar\eps_i,
\]
with $\sum_i\bar\eps_i<\infty$.  For $j>k$ choose $\ell$ as in \Cref{lem:cross}.  Using
$t_i\le\bar t$, \Cref{lem:cross}, and
$\norm{x_\ell-x^\dagger}\le R_*$ gives
\[
 \left|\ip{x_\ell-x_k}{x_\ell-x^\dagger}\right|
 \le C\sum_{i=k}^{\ell-1}t_i\norm{r_i}^2+R_*\sum_{i=k}^\infty\bar\eps_i\to0.
\]
Since $\norm{x_i-x^\dagger}^2$ converges, the identity
\[
 \norm{x_\ell-x_k}^2
 =\norm{x_k-x^\dagger}^2-\norm{x_\ell-x^\dagger}^2
 +2\ip{x_\ell-x_k}{x_\ell-x^\dagger}
\]
shows $\norm{x_\ell-x_k}\to0$ uniformly in $j>k$.  For the second segment,
the same increment expansion, now from $\ell$ to $j$, gives
\[
 \left|\ip{x_j-x_\ell}{x_\ell-x^\dagger}\right|
 \le C\sum_{i=\ell}^{j-1}t_i\norm{r_i}^2
      +R_*\sum_{i=\ell}^\infty\bar\eps_i\longrightarrow0
\]
uniformly in $j>k$ as $k\to\infty$.  Since
\[
 \norm{x_j-x_\ell}^2
 =\norm{x_j-x^\dagger}^2-\norm{x_\ell-x^\dagger}^2
 -2\ip{x_j-x_\ell}{x_\ell-x^\dagger},
\]
convergence of $\norm{x_i-x^\dagger}^2$ implies
$\norm{x_j-x_\ell}\to0$ uniformly.  The triangle inequality now gives
$\sup_{j>k}\norm{x_j-x_k}\to0$.  Continuity of $F$ and $r_k\to0$
identify the limit as a solution.
\end{proof}

\begin{lemma}[Quasi-Fej\'er estimate relative to a local solution]\label{lem:qf-anysolution}
Let \Cref{ass:tcc,ass:margin} hold, and let $\bar x\in\overline B_{R_*}(x^\dagger)$ satisfy $F(\bar x)=y$.  Then, for all sufficiently small $\delta$, there is a nonnegative summable sequence $\{\xi_k(\bar x)\}$, independent of $\delta$, such that every pre-discrepancy iterate satisfies
\begin{equation}\label{eq:qf-anysolution}
 \norm{x_{k+1}^\delta-\bar x}^2
 \le \norm{x_k^\delta-\bar x}^2
 -\kappa_\tau\norm{r_k^\delta}^2+\xi_k(\bar x).
\end{equation}
In particular, one may take
\[
 C_{\bar x}:=\norm{x_0-\bar x}+\Omega\bar\sigma_0,
 \qquad
 \xi_k(\bar x):=2C_{\bar x}\bar\eps_k+\bar\eps_k^2.
\]
\end{lemma}
\begin{proof}
The localization already established in \Cref{lem:qf} keeps $x_k^\delta$ in $\overline B_{R_*}(x^\dagger)$, while $\bar x$ is a solution in the same TCC neighborhood.  Repeating the proof of \Cref{lem:basicinner} with $\bar x$ in place of $x^\dagger$ therefore gives
\[
 \norm{z_{k+1}^\delta-\bar x}^2
 \le \norm{x_k^\delta-\bar x}^2
 -\kappa_\tau\norm{r_k^\delta}^2.
\]
Since $\norm{e_k^\delta}\le\bar\eps_k$, induction yields
$\norm{x_k^\delta-\bar x}\le C_{\bar x}$ and hence also
$\norm{z_{k+1}^\delta-\bar x}\le C_{\bar x}$.  Expanding
$x_{k+1}^\delta=z_{k+1}^\delta+e_k^\delta$ gives \eqref{eq:qf-anysolution}.  Summability follows from $\sum_k\omega_k<\infty$.
\end{proof}

\subsection{Stability at fixed iteration indices}

\begin{lemma}[Stability at fixed iteration indices]\label{lem:AAcontinuity}
Let \Cref{ass:tcc,ass:margin} hold, and let $y^{\delta_n}\to y$ with
$\norm{y^{\delta_n}-y}\le\delta_n$.  For every fixed $K\ge0$, the constant extensions of the noisy and exact trajectories satisfy
\begin{equation}\label{eq:fixedstability}
 x_k^{\delta_n}\longrightarrow x_k,
 \qquad 0\le k\le K.
\end{equation}
\end{lemma}

\begin{proof}
We argue by induction.  The assertion is immediate for $k=0$.  Suppose it holds through index $k$.  If $F(x_k)\ne y$, the noisy trajectories cannot stop before index $k$ for infinitely many $n$: otherwise a subsequence would stop at one common earlier index, and the induction hypothesis together with \eqref{eq:discrepancy} would imply that the corresponding exact iterate is already a solution.  Moreover, \Cref{lem:basicinner} gives $g_k\ne0$.  Hence the residual, the certified step length, and the anchor are continuous at the exact history.  Since $f_0\ne0$, we also have
$\lambda_{\rm AA}^0=\eps_{\rm AA}\norm{f_0}^2>0$ and
\[
 \gamma_k^\delta=
 \bigl((\Delta F_k^\delta)^*\Delta F_k^\delta
       +\lambda_{\rm AA}^\delta I\bigr)^{-1}
 (\Delta F_k^\delta)^*f_k^\delta
\]
depends continuously on the finite history and the data.  The raw AA proposal and the radial projection in \eqref{eq:correction} are therefore continuous, proving $x_{k+1}^{\delta_n}\to x_{k+1}$.

If $F(x_k)=y$, the exact trajectory is constant from index $k$ onward.  For any noisy trajectory that has not already stopped,
$\norm{z_{k+1}^{\delta_n}-x_k^{\delta_n}}
\le\bar t\norm{g_k^{\delta_n}}\to0$ and
$\norm{e_k^{\delta_n}}\le\chi\sigma_k^{\delta_n}\to0$.
Thus $x_{k+1}^{\delta_n}\to x_k=x_{k+1}$ without requiring continuity of the raw proposal at the zero-residual history.  This closes the induction.
\end{proof}

\subsection{Finite discrepancy stopping}

\begin{theorem}[Finite discrepancy stopping]\label{thm:finite}
Under \Cref{ass:tcc,ass:margin}, for every $0<\delta\le\delta_0$ the discrepancy index
\[
 k_\delta=\min\{k:\norm{F(x_k^\delta)-y^\delta}\le\tau\delta\}
\]
is finite, and
\begin{equation}\label{eq:crudestop}
 k_\delta\le 1+\frac{d_0^2+\theta}{\kappa_\tau\tau^2\delta^2}.
\end{equation}
\end{theorem}
\begin{proof}
If the first $N$ iterates all precede discrepancy stopping, then \eqref{eq:qf} and $\norm{r_k^\delta}>\tau\delta$ give
\[
 0\le d_0^2-\kappa_\tau N\tau^2\delta^2+\theta.
\]
Solving for $N$ proves the claim.
\end{proof}

\subsection{Iterative regularization}

\begin{theorem}[Strong iterative regularization without local uniqueness]\label{thm:reg}
Let \Cref{ass:tcc,ass:margin} hold, and let $x_*$ be the strong limit of the exact-data trajectory from \Cref{thm:strong}.  Then the discrepancy-stopped outputs of Algorithm~\ref{alg:rglAA} satisfy
\[
 x_{k_\delta}^\delta\to x_*\qquad\text{strongly as }\delta\to0.
\]
Thus the noisy method regularizes toward the same solution selected by the exact-data dynamics.
\end{theorem}
\begin{proof}
Let $\delta_n\to0$.  If a subsequence of $k_{\delta_n}$ is bounded, pass to a further subsequence on which $k_{\delta_n}=K$.  By \Cref{lem:AAcontinuity}, $x_K^{\delta_n}\to x_K$, while the discrepancy principle gives $F(x_K)=y$.  For exact data the algorithm stops once a zero residual is reached and is then extended constantly, so $x_K=x_*$.

It remains to consider a subsequence for which $k_{\delta_n}\to\infty$.  Since $x_*\in\overline B_{R_*}(x^\dagger)$ and $F(x_*)=y$, \Cref{lem:qf-anysolution} applies with $\bar x=x_*$.  For every fixed $K$ and all sufficiently large $n$, summing \eqref{eq:qf-anysolution} from $K$ to $k_{\delta_n}-1$ and dropping the residual decrease gives
\[
 \norm{x_{k_{\delta_n}}^{\delta_n}-x_*}^2
 \le \norm{x_K^{\delta_n}-x_*}^2
 +\sum_{j=K}^\infty\xi_j(x_*).
\]
First let $n\to\infty$ and use \Cref{lem:AAcontinuity}; then let $K\to\infty$ and use $x_K\to x_*$ from \Cref{thm:strong} together with the summability of $\xi_j(x_*)$.  Every sequence $\delta_n\to0$ therefore converges to $x_*$.
\end{proof}

\begin{corollary}[Identification with a prescribed solution]\label{cor:unique}
Under the assumptions of \Cref{thm:reg}, if $x^\dagger$ is the only solution of $F(x)=y$ in $\overline B_{R_*}(x^\dagger)$, then $x_*=x^\dagger$ and hence $x_{k_\delta}^\delta\to x^\dagger$ strongly.
\end{corollary}

\section{Benchmark convergence rate}
This section treats the benchmark source condition corresponding to the classical Hilbert-space exponent $\nu=1/2$.  All assumptions used in the rate analysis are stated before the first auxiliary result.

\subsection{Additional assumptions and proof strategy}

Assumptions~\ref{ass:tcc} and~\ref{ass:margin} remain in force, the memory $m$ is fixed, and the AA proposal uses the fixed-reference ridge \eqref{eq:ridge}.  Set
\begin{equation}\label{eq:Tdef}
 T:=F'(x^\dagger).
\end{equation}

\begin{assumption}[Range invariance and benchmark source condition]\label{ass:rate}
There exist bounded operators $Q_x:\Y\to\Y$ on $\overline B_{R_*}(x^\dagger)$ such that
\begin{equation}\label{eq:rangeinv}
 F'(x)=Q_xT,
 \qquad
 \norm{I-Q_x}\le\kappa_0\norm{x-x^\dagger}.
\end{equation}
The constants satisfy $\kappa_0R_*\le q<1$.  Moreover,
\begin{equation}\label{eq:benchmarksource}
 x^\dagger-x_0=T^*\lambda^\dagger
 \qquad\text{for some }\lambda^\dagger\in\Y.
\end{equation}
\end{assumption}

\begin{assumption}[Rate-compatible outer schedule]\label{ass:ratetail}
In addition to \eqref{eq:omegasum}, the correction budget satisfies
\begin{equation}\label{eq:tailomega}
 C_\omega:=\sup_{n\ge1}n\sum_{k=n}^\infty\omega_k<\infty.
\end{equation}
\end{assumption}

The factorization in \eqref{eq:rangeinv}, including the Lipschitz closeness of $Q_x$ to the identity, is exactly the range-invariance hypothesis used in Assumption~3 of \citet{Jin2025}, specialized here to Hilbert spaces and written with $T=F'(x^\dagger)$.  The explicit bound $q<1$ follows after restricting the local ball so that $\kappa_0R_*<1$.  Equation~\eqref{eq:benchmarksource} is the corresponding benchmark Hilbert-space source condition.  No local uniqueness assumption is imposed in the rate analysis; the source and range-invariance conditions select the reference solution $x^\dagger$ directly.

The proof has three steps.  First, range invariance keeps every iterate in $x_0+\operatorname{Ran}(T^*)$ and converts the accepted AA correction into a source-space perturbation.  Second, the ridge and finite memory bound this perturbation by the residual energy from Section~3, which yields source-element stability up to $O(\delta^{-1})$ iterations.  Third, the negative term in the source-energy estimate and the tail condition \eqref{eq:tailomega} give the $O(\delta^{1/2})$ discrepancy-principle rate.

\subsection{Source representation and the Anderson dual perturbation}
\begin{lemma}[Source representation]\label{lem:sourcerep}
Under \Cref{ass:rate}, there are elements $\lambda_k^\delta\in\Y$ such that
\begin{equation}\label{eq:sourcerep}
 x_k^\delta=x_0+T^*\lambda_k^\delta.
\end{equation}
If
\[
 \nu_k^\delta:=\lambda_k^\delta-\mu Q_k^*r_k^\delta,
 \qquad Q_k:=Q_{x_k^\delta},
\]
then the AA proposal has source element
\begin{equation}\label{eq:aasource}
 \widehat\lambda_{k+1}^{\rm AA,\delta}
 =\nu_k^\delta-\Delta N_k^\delta\gamma_k^\delta,
\end{equation}
where $\Delta N_k^\delta$ contains the recent differences of $\nu_j^\delta$.  The committed source recursion is
\begin{equation}\label{eq:sourcerec}
 \lambda_{k+1}^\delta
 =\lambda_k^\delta-\bar t_k^\delta Q_k^*r_k^\delta+q_k^\delta,
\end{equation}
with
\begin{equation}\label{eq:tbarq}
 \bar t_k^\delta=(1-\beta_k^\delta)t_k^\delta+\beta_k^\delta\mu,
 \qquad
 q_k^\delta=-\beta_k^\delta\Delta N_k^\delta\gamma_k^\delta.
\end{equation}
Before discrepancy stopping,
\begin{equation}\label{eq:tbarbounds}
 0<h_-\le\bar t_k^\delta\le h_+<\infty
\end{equation}
with
\[
 h_-:=\min\{t_*,\mu\},\qquad h_+:=\max\{\bar t,\mu\},
\]
independent of $k$ and sufficiently small $\delta$.
\end{lemma}
\begin{proof}
Set $\lambda_0^\delta=0$.  Because $F'(x)^*=T^*Q_x^*$, both the
certified anchor and the raw map preserve the affine space
$x_0+\operatorname{Ran}(T^*)$.  If the stored raw-map source elements are
$\nu_j^\delta$, taking their same affine differences as in
\eqref{eq:aa-cand} gives \eqref{eq:aasource}.  Substituting the anchor and
AA source representations into
$x_{k+1}^\delta=(1-\beta_k^\delta)z_{k+1}^\delta+
\beta_k^\delta \widehat{x}_{k+1}^{\rm AA,\delta}$ and collecting source elements
gives \eqref{eq:sourcerec}--\eqref{eq:tbarq}.  Finally,
$0\le\beta_k^\delta\le1$, $t_*\le t_k^\delta\le\bar t$, and $\mu>0$
give the stated bounds.
\end{proof}

For the next estimates let
\begin{equation}\label{eq:Sk}
 S_k^\delta:=\sum_{j=\max\{0,k-m\}}^k\norm{r_j^\delta}^2,
 \qquad
 A_k^\delta:=\max_{0\le j\le k}\norm{\lambda_j^\delta-\lambda^\dagger}.
\end{equation}
By finite memory and \eqref{eq:resenergy},
\begin{equation}\label{eq:Ssum}
 \sum_{k<N}S_k^\delta\le (m+1)C_r
 \qquad\text{for }N\le k_\delta.
\end{equation}

\begin{lemma}[AA coefficient and dual-correction bounds]\label{lem:gammabound}
For all sufficiently small $\delta$, the fixed-reference ridge satisfies $\lambda_{\rm AA}^\delta\ge\lambda_*>0$.  Moreover,
\begin{equation}\label{eq:gammabound}
 \norm{\gamma_k^\delta}\le C S_k^\delta,
\end{equation}
and
\begin{equation}\label{eq:qbound}
 \norm{q_k^\delta}\le C S_k^\delta(A_k^\delta+1).
\end{equation}
\end{lemma}
\begin{proof}
The regularized AA normal equations give
\[
 \norm{\gamma_k^\delta}
 \le \frac{\norm{\Delta F_k^\delta}\,\norm{f_k^\delta}}{\lambda_{\rm AA}^\delta}.
\]
Because $f_j^\delta=-\mu F'(x_j^\delta)^*r_j^\delta$ and the derivative is uniformly bounded, finite memory gives
\[
 \norm{\Delta F_k^\delta}\le C(S_k^\delta)^{1/2},
 \qquad \norm{f_k^\delta}\le C(S_k^\delta)^{1/2},
\]
proving \eqref{eq:gammabound}.  The source candidates $\nu_j^\delta=\lambda_j^\delta-\mu Q_j^*r_j^\delta$ satisfy $\norm{\nu_j^\delta}\le C(A_k^\delta+1)$ on the local ball.  Hence the finite difference matrix $\Delta N_k^\delta$ has norm at most $C(A_k^\delta+1)$, and \eqref{eq:qbound} follows from \eqref{eq:tbarq}.
\end{proof}

\subsection{Source-element stability}
The rate argument uses two distinct consequences of range invariance.  We isolate them first, because they are the point at which the nonlinear residual is transferred to the fixed linearized range space.

\begin{lemma}[Range-invariance linearization and residual transfer]\label{lem:range-transfer}
Let \Cref{ass:rate} hold and set $T=F'(x^\dagger)$.  For $x\in\overline B_{R_*}(x^\dagger)$ and $e=x-x^\dagger$,
\begin{equation}\label{eq:lin-defect}
 \norm{F(x)-y-Te}\le \frac{\kappa_0}{2}\norm e\,\norm{Te},
 \qquad
 \norm{Te}\le \frac{1}{1-q}\norm{F(x)-y}.
\end{equation}
If $r^\delta=F(x)-y^\delta$ is pre-discrepancy, then
\begin{equation}\label{eq:res-transfer}
 Q_x^*r^\delta=Te+d,
 \qquad
 \norm d\le C\norm e\,\norm{r^\delta}+C\delta,
 \qquad
 \norm{Q_x^*r^\delta}\le C\norm{r^\delta}.
\end{equation}
\end{lemma}
\begin{proof}
The fundamental theorem of calculus and \eqref{eq:rangeinv} give
\[
 F(x)-y-Te=\int_0^1 (Q_{x^\dagger+se}-I)Te\,ds,
\]
which proves the first estimate in \eqref{eq:lin-defect}.  Since $\kappa_0\norm e\le q<1$, the second estimate follows by the triangle inequality.  Next write
\[
 Q_x^*r^\delta-Te=(Q_x^*-I)r^\delta+\bigl(F(x)-y-Te\bigr)-(y^\delta-y).
\]
The first term is bounded by $\kappa_0\norm e\,\norm{r^\delta}$.  The nonlinear defect is controlled by \eqref{eq:lin-defect}, the second bound in \eqref{eq:lin-defect}, and
$\norm{F(x)-y}\le\norm{r^\delta}+\delta$.  Before discrepancy stopping, $\delta<\norm{r^\delta}/\tau$, so all lower-order $\delta$ terms may be absorbed into the stated estimate for $d$.  The final bound in \eqref{eq:res-transfer} follows from $Q_x^*r^\delta=Te+d$ and the preceding estimates.
\end{proof}

The following one-step estimate is the Hilbert-space analogue of the key source-energy inequality in \citet{Jin2025}.

\begin{lemma}[Uniform source one-step estimate]\label{lem:sourcestep}
Under \Cref{ass:rate}, there is $c_e>0$ such that for every $h\in[h_-,h_+]$ and every pre-discrepancy iterate,
\begin{align}\label{eq:sourcestep}
 &\norm{\lambda_k^\delta-hQ_k^*r_k^\delta-\lambda^\dagger}^2\notag\\
 &\quad\le
 (1+C\norm{r_k^\delta}^2)\norm{\lambda_k^\delta-\lambda^\dagger}^2
 +C\norm{r_k^\delta}^2
 +C\delta\norm{\lambda_k^\delta-\lambda^\dagger}
 -c_e\norm{x_k^\delta-x^\dagger}^2.
\end{align}
\end{lemma}
\begin{proof}
Let $e=x_k^\delta-x^\dagger=T^*(\lambda_k^\delta-\lambda^\dagger)$ and set $a=\lambda_k^\delta-\lambda^\dagger$.  By \Cref{lem:range-transfer},
\[
 Q_k^*r_k^\delta=Te+d_k,
 \qquad
 \norm{d_k}\le C\norm e\,\norm{r_k^\delta}+C\delta,
 \qquad
 \norm{Q_k^*r_k^\delta}\le C\norm{r_k^\delta}.
\]
Expanding the square gives the leading term
\[
 -2h\ip{a}{Te}=-2h\norm e^2.
\]
The cross term involving $d_k$ is bounded by Young's inequality,
\[
 C\norm a\,\norm e\,\norm{r_k^\delta}
 \le \frac{c_e}{2}\norm e^2+C\norm{r_k^\delta}^2\norm a^2,
\]
and the noise term contributes $C\delta\norm a$.  The last estimate in \eqref{eq:res-transfer} bounds the quadratic step term by $C\norm{r_k^\delta}^2$.  The uniform bounds
$h_-\le h\le h_+$ leave a fixed negative multiple of $\norm e^2$,
proving \eqref{eq:sourcestep}.
\end{proof}

\begin{theorem}[Source-element stability up to $O(\delta^{-1})$ iterations]\label{thm:sourcestability}
Let \Cref{ass:tcc,ass:margin,ass:rate} hold and let the AA memory be fixed.  For every fixed $c_*>0$ there exists $C>0$ such that
\begin{equation}\label{eq:sourcestable}
 \norm{\lambda_k^\delta-\lambda^\dagger}\le C,
 \qquad
 0\le k\le\min\left\{k_\delta,\left\lfloor\frac{c_*}{\delta}\right\rfloor\right\},
\end{equation}
for all sufficiently small $\delta$.
\end{theorem}
\begin{proof}
Apply \Cref{lem:sourcestep} with $h=\bar t_k^\delta$ and write
$b_k^\delta=\lambda_k^\delta-\lambda^\dagger-
\bar t_k^\delta Q_k^*r_k^\delta$.  On the local ball,
$\norm{b_k^\delta}\le C(A_k^\delta+1)$.  Thus \eqref{eq:qbound} and
$2\ip{b_k^\delta}{q_k^\delta}+\norm{q_k^\delta}^2$ contribute at most
$CS_k^\delta((A_k^\delta)^2+1)$; here we used the uniform boundedness of
$S_k^\delta$ following from \eqref{eq:resenergy}.  Dropping the favorable
$-c_e\norm{x_k^\delta-x^\dagger}^2$ term gives
\[
 (A_{k+1}^\delta)^2
 \le (1+CS_k^\delta)(A_k^\delta)^2+CS_k^\delta+C\delta A_k^\delta.
\]
The discrete Grönwall inequality and \eqref{eq:Ssum} yield
\[
 (A_n^\delta)^2\le C_0+C_1\delta\sum_{j=0}^{n-1}A_j^\delta
 \le C_0+C_1\delta n A_n^\delta.
\]
Solving the scalar quadratic inequality gives
\[
 A_n^\delta\le \sqrt{C_0}+C_1\delta n.
\]
Hence $A_n^\delta$ is uniformly bounded whenever $n\le c_*/\delta$.
\end{proof}

\subsection{A priori and discrepancy-principle benchmark rate}
Once \Cref{thm:sourcestability} is available, the favorable term in \eqref{eq:sourcestep} can be retained.

\begin{lemma}[Accumulated primal error]\label{lem:sumprimal}
For $N\le\min\{k_\delta,\lfloor c_*/\delta\rfloor\}$,
\begin{equation}\label{eq:sumprimal}
 \sum_{k=0}^{N-1}\norm{x_k^\delta-x^\dagger}^2\le C.
\end{equation}
\end{lemma}
\begin{proof}
By \Cref{thm:sourcestability}, $A_k^\delta$ is uniformly bounded in the stated range.  Hence \eqref{eq:qbound} reduces to $\norm{q_k^\delta}\le CS_k^\delta$.  With $b_k^\delta$ as above,
\[
 2|\ip{b_k^\delta}{q_k^\delta}|+\norm{q_k^\delta}^2
 \le C S_k^\delta,
\]
again because $(S_k^\delta)^2\le C S_k^\delta$.  Keeping the negative
term in \eqref{eq:sourcestep} gives
\[
 \norm{\lambda_{k+1}^\delta-\lambda^\dagger}^2
 +c\norm{x_k^\delta-x^\dagger}^2
 \le
 \norm{\lambda_k^\delta-\lambda^\dagger}^2
 +CS_k^\delta+C\delta.
\]
Sum from $0$ to $N-1$ and use \eqref{eq:Ssum} and $\delta N\le c_*$.
\end{proof}

\begin{lemma}[Rate transfer from averaging]\label{lem:averaging}
Under \Cref{ass:ratetail}, if $N\le k_\delta$ and $N\le c_*/\delta$, then
\begin{equation}\label{eq:ENrate}
 \norm{x_N^\delta-x^\dagger}^2\le \frac{C}{N+1}.
\end{equation}
\end{lemma}
\begin{proof}
The result is immediate after increasing $C$ when $N<2$.  For $N\ge2$,
write $E_k=\norm{x_k^\delta-x^\dagger}^2$.  Dropping the negative term in the quasi-Fej\'er inequality gives, for $L\le k<N$,
\[
 E_N\le E_k+\sum_{j=k}^{N-1}\xi_j.
\]
Take $L=\lfloor N/2\rfloor$ and average over $k=L,\ldots,N-1$.  By \Cref{lem:sumprimal}, the average of the $E_k$ terms is $O(N^{-1})$.  Moreover, boundedness of the summable sequence $\{\omega_k\}$ gives $\omega_k^2\le C\omega_k$; hence \eqref{eq:tailomega} and \eqref{eq:bareps} imply
\[
 \sum_{j=L}^\infty\xi_j\le \frac{C}{L+1}\le\frac{C}{N+1}.
\]
This proves \eqref{eq:ENrate}.
\end{proof}

\begin{theorem}[Benchmark discrepancy-principle rate]\label{thm:sqrtrate}
Let \Cref{ass:tcc,ass:margin,ass:rate,ass:ratetail} hold.  Then
\begin{equation}\label{eq:sqrtrate}
 \norm{x_{k_\delta}^\delta-x^\dagger}\le C\delta^{1/2}
\end{equation}
for all sufficiently small $\delta$.
\end{theorem}
\begin{proof}
Fix $c_*>0$ and, for sufficiently small $\delta$, set
$N_\delta=\lfloor c_*/\delta\rfloor\ge2$.  If $k_\delta\le N_\delta$, then \Cref{thm:sourcestability} gives a uniform bound on $\norm{\lambda_{k_\delta}^\delta-\lambda^\dagger}$.  By \eqref{eq:sourcerep} and \eqref{eq:benchmarksource},
\[
 \norm{x_{k_\delta}^\delta-x^\dagger}^2
 =\ip{\lambda_{k_\delta}^\delta-\lambda^\dagger}{T(x_{k_\delta}^\delta-x^\dagger)}.
\]
For $e=x-x^\dagger$, the direct integral estimate used in
\Cref{lem:sourcestep} yields
\[
 \norm{F(x)-y-Te}\le q\norm{Te},\qquad q<1.
\]
Consequently,
\[
 \norm{T(x_{k_\delta}^\delta-x^\dagger)}
 \le C\norm{F(x_{k_\delta}^\delta)-y}
 \le C(\tau+1)\delta.
\]
Hence the squared error is $O(\delta)$.

If $k_\delta>N_\delta$, then all first $N_\delta$ iterates are pre-discrepancy.  By \Cref{lem:averaging},
\[
 \norm{x_{N_\delta}^\delta-x^\dagger}^2\le C/N_\delta\le C\delta.
\]
Using the quasi-Fej\'er inequality from $N_\delta$ to $k_\delta$ and the tail condition once more,
\[
 \norm{x_{k_\delta}^\delta-x^\dagger}^2
 \le C\delta+\sum_{j=N_\delta}^\infty\xi_j
 \le C\delta.
\]
Taking square roots proves \eqref{eq:sqrtrate}.
\end{proof}
%

\section{Numerical experiments}\label{sec:numerics}
\begingroup
\normalsize
\setlength{\textfloatsep}{7pt plus 2pt minus 2pt}
\setlength{\floatsep}{7pt plus 2pt minus 2pt}
\setlength{\intextsep}{7pt plus 2pt minus 2pt}
\captionsetup{font=small,skip=3pt}
\renewcommand{\arraystretch}{0.94}

We present numerical results for \RG on nonlinear inverse problems. The diagonal example examines the benchmark convergence rate and compares residual and coefficient histories. The PDE examples compare reconstruction accuracy and computational work, and examine the effects of parameter choices, observation regions, and mesh refinement. To reduce the influence of random noise fluctuations, all reported means are computed over 20 noise realizations, shared across methods within each test case. For exact discrete data $y$, noisy data are generated by
\[
 y^\delta=y+\delta\frac{Z}{\|Z\|_{Y,h}},
 \qquad \delta=\delta_{\mathrm{rel}}\|y\|_{Y,h},
\]
where the components of $Z$ are independent standard normal random variables, $\delta_{\mathrm{rel}}$ is the relative noise level, and $\|\cdot\|_{Y,h}$ denotes the discrete data norm.

Unless otherwise stated, we set
\[
 \tau=1.5,\qquad \bar\gamma=1,\qquad m=5,\qquad
 \eps_{\rm AA}=10^{-10},\qquad \Omega=10,\qquad s=2,\qquad \chi=5.
\]
We use $\eta=0.12$, except that $\eta=0.15$ in the two-dimensional potential tests. We set the step-size cap to $\bar t=10^{12}$; this bound is never reached in the experiments. For comparison, the unsafeguarded baseline, AA--Landweber, uses zero ridge and LSQR, whereas \RG uses the fixed-reference ridge with $\eps_{\rm AA}=10^{-10}$ and normal equations. For each problem family, the proposal step $\mu$ is selected on independent pilot data at $\delta_{\mathrm{rel}}=10^{-3}$, then fixed and shared by the two AA methods.

To distinguish different aspects of computational cost in the tables, ``PDE solves'' counts forward, adjoint, and, for IRGNM, inner linearized solves, whereas ``factorizations'' counts sparse matrix factorizations that may be reused across solves. The ``wall time'' column reports the implementation-dependent running time. Runs used Python 3.12.14, NumPy 2.3.5, SciPy 1.17.0, pandas 3.0.5, and Matplotlib 3.11.2 on x86--64 macOS.

\subsection{Parameter sensitivity on a one-dimensional elliptic problem}\label{sec:parameter-sensitivity}
We first examine the sensitivity of \RG to the memory depth $m$, decay exponent $s$, ridge parameter $\eps_{\rm AA}$, total correction budget $\Omega$, and relative safeguard parameter $\chi$. Consider recovering $q$ from full internal observations $F(q)=u(q)$, where
\[
 -u''+qu=10\quad\text{in }(0,1),\qquad u(0)=u(1)=0.
\]
For a perturbation $\varphi$, the linearized state $v=F'(q)\varphi$ and the adjoint state $p$ satisfy
\[
 (-\partial_{xx}+q)v=-\varphi u,
 \qquad
 (-\partial_{xx}+q)p=r,
 \qquad
 F'(q)^*r=-up.
\]
We use 1024 interior grid points, $q_0\equiv0.8$, and the rough exact potential
\begin{equation}\label{eq:1d-rough-truth}
 q_{\mathrm{rough}}^\dagger(x)=
 \begin{cases}
  3.8, & 0.20<x<0.38,\\
  2.5, & 0.58<x<0.70,\\
  4.8, & 0.80<x<0.88,\\
  0.6, & \text{otherwise}.
 \end{cases}
\end{equation}
The truth is shown in the right panel of \Cref{fig:1d-truths}. At $\delta_{\mathrm{rel}}=10^{-4}$, each parameter is varied separately while the others remain fixed at their default values.

Table~\ref{tab:parameter-sensitivity} reports the mean PDE-solve count, its sample standard deviation (SD), and the mean reconstruction error. The column $\gamma_{\max}$ reports the average maximum AA coefficient norm, obtained by averaging $\max_{1\le k<k_\delta}\|\gamma_k^\delta\|_2$ over noise realizations. The column corr. reports the average cumulative accepted correction, obtained by averaging $\sum_{k=0}^{k_\delta-1}\|e_k^\delta\|_h$ over noise realizations. Here, $e_k^\delta$ is the accepted correction defined in Section~2, and the one-dimensional discrete $L^2$ norm is
\[
 \|v\|_h=\left(h\sum_i|v_i|^2\right)^{1/2},\qquad h=1/1025,
\]
where $h$ is the grid spacing.

\begin{table}[H]
\centering
\renewcommand{\arraystretch}{1.00}
\caption{One-factor-at-a-time sensitivity over 20 common-noise realizations.}
\label{tab:parameter-sensitivity}
\footnotesize
\setlength{\tabcolsep}{4pt}
\begin{tabular}{llrrrrr}
\toprule
Parameter & value & PDE solves & SD & error & $\gamma_{\max}$ & corr.\\
\midrule
$m$ & $1$ & 987.4 & 158.7 & 0.57375 & 1.0 & 1.866\\
 & $2$ & 1006.6 & 196.4 & 0.57340 & 20.0 & 1.688\\
 & $5$ & 1000.6 & 195.8 & 0.57210 & 35.5 & 1.720\\
 & $10$ & 1135.8 & 229.9 & 0.57695 & 34.9 & 1.722\\
\midrule
$s$ & $1.5$ & 697.1 & 113.3 & 0.57475 & 34.9 & 2.986\\
 & $2$ & 1000.6 & 195.8 & 0.57210 & 35.5 & 1.720\\
 & $2.5$ & 1411.9 & 355.5 & 0.57541 & 35.4 & 1.122\\
 & $3$ & 1369.6 & 298.0 & 0.57715 & 35.9 & 0.825\\
\midrule
$\eps_{\rm AA}$ & $10^{-16}$ & 1275.3 & 244.9 & 0.57373 & 2352.4 & 1.723\\
 & $10^{-14}$ & 1222.5 & 224.0 & 0.57575 & 741.1 & 1.722\\
 & $10^{-12}$ & 1003.7 & 140.4 & 0.57513 & 195.3 & 1.721\\
 & $10^{-10}$ & 1000.6 & 195.8 & 0.57210 & 35.5 & 1.720\\
 & $10^{-8}$ & 1072.0 & 254.6 & 0.57324 & 26.4 & 1.710\\
\midrule
$\Omega$ & $1$ & 1346.9 & 324.8 & 0.57630 & 28.1 & 0.541\\
 & $5$ & 1143.1 & 187.5 & 0.57444 & 55.7 & 1.197\\
 & $10$ & 1000.6 & 195.8 & 0.57210 & 35.5 & 1.720\\
 & $15$ & 954.6 & 215.2 & 0.57377 & 35.3 & 2.184\\
 & $20$ & 846.1 & 88.6 & 0.57247 & 34.9 & 2.567\\
 & $30$ & 761.3 & 122.2 & 0.57336 & 35.0 & 3.188\\
 & $40$ & 667.3 & 135.9 & 0.57327 & 34.9 & 3.851\\
 & $50$ & 620.6 & 114.3 & 0.57520 & 34.9 & 4.363\\
\midrule
$\chi$ & $1$ & 1050.9 & 206.7 & 0.57629 & 35.4 & 1.698\\
 & $2$ & 1000.6 & 195.8 & 0.57210 & 35.5 & 1.720\\
 & $5$ & 1000.6 & 195.8 & 0.57210 & 35.5 & 1.720\\
 & $10$ & 1000.6 & 195.8 & 0.57210 & 35.5 & 1.720\\
\bottomrule
\end{tabular}
\end{table}

Table~\ref{tab:parameter-sensitivity} shows that the reconstruction errors vary little over the tested parameter values, whereas the computational work, Anderson coefficient norms, and cumulative accepted corrections are more sensitive to these choices. Increasing the memory depth does not consistently reduce the PDE-solve count. Stronger ridge regularization reduces the coefficient norms, but does not lead to a monotonic reduction in computational work. Increasing the total correction budget allows larger cumulative accepted corrections and reduces the PDE-solve count. Among the tested rate-compatible decay exponents, $s=2$ requires the fewest PDE solves. The reported summaries are identical for the tested values $\chi=2,5,10$.

For the subsequent comparisons, we retain $(m,s,\eps_{\rm AA},\Omega,\chi)=(5,2,10^{-10},10,5)$ as the common default unless otherwise stated. Here, $m=5$ and $\Omega=10$ provide intermediate settings within the tested ranges, while $\eps_{\rm AA}=10^{-10}$ combines controlled coefficient growth with low computational work. The choice $s=2$ also satisfies the tail condition required for the convergence-rate analysis.

\subsection{One-dimensional elliptic potential comparison}
Using the model and initial guess from Section~\ref{sec:parameter-sensitivity}, we compare three methods on 1024 interior grid points: safe Landweber, AA--Landweber, and \RG. The rough truth is given by \eqref{eq:1d-rough-truth}, and the smooth truth is
\begin{equation}\label{eq:1d-smooth-truth}
 q_{\mathrm{smooth}}^\dagger(x)
 =0.7+2.8\exp\!\left[-\left(\frac{x-0.33}{0.12}\right)^2\right]
 +1.8\exp\!\left[-\left(\frac{x-0.72}{0.09}\right)^2\right].
\end{equation}
Both truths are shown in \Cref{fig:1d-truths} and tested at $\delta_{\mathrm{rel}}=10^{-3}$ and $10^{-4}$.

\begin{figure}[H]
\centering
\includegraphics[width=.88\textwidth]{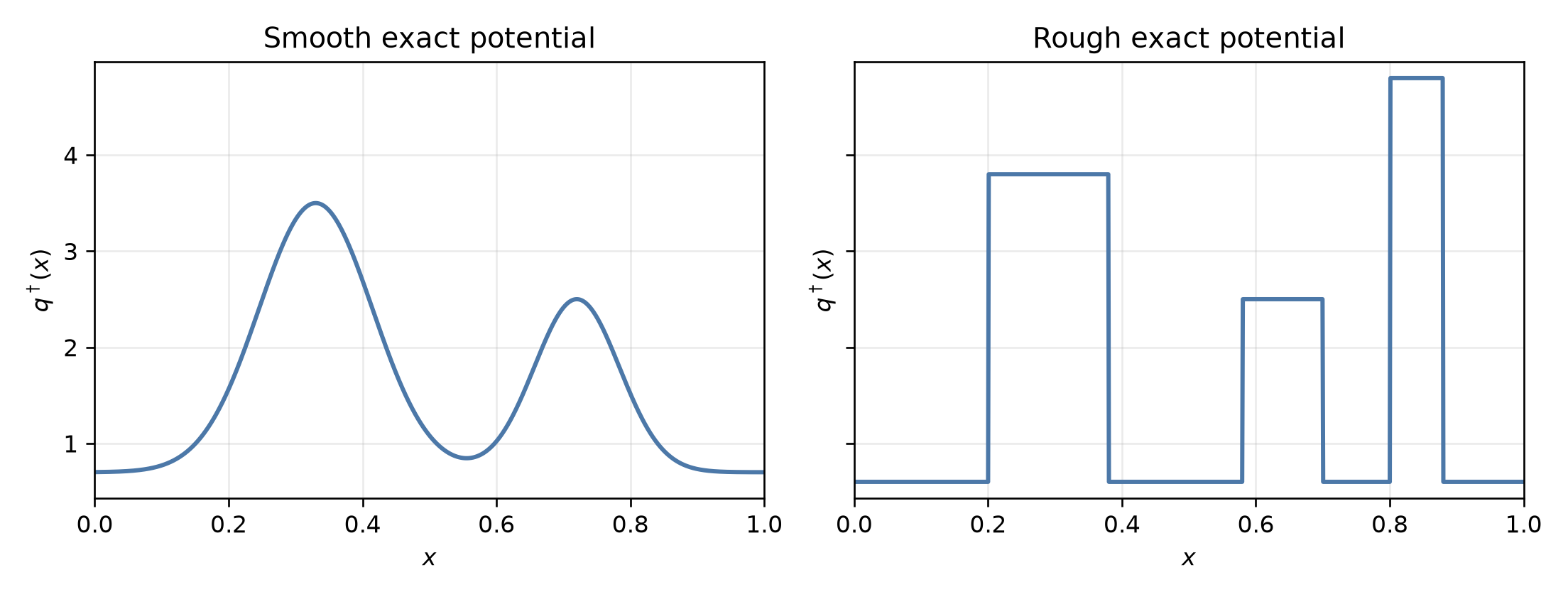}
\caption{One-dimensional exact potentials: smooth (left) and rough (right).}
\label{fig:1d-truths}
\end{figure}

\begin{table}[H]
\centering
\renewcommand{\arraystretch}{1.00}
\caption{Mean results for one-dimensional potential inversion at two relative-noise levels.}
\label{tab:1d-final}
\small
\begin{tabular}{lllrr}
\toprule
Truth & $\delta_{\mathrm{rel}}$ & Method & PDE solves & error\\
\midrule
smooth & $10^{-3}$ & Safe Landweber & 86.5 & 0.29725\\
 & & AA--Landweber & 31.5 & 0.16656\\
 & & \RG & 73.0 & 0.31211\\
\midrule
smooth & $10^{-4}$ & Safe Landweber & 301.3 & 0.06904\\
 & & AA--Landweber & 55.3 & 0.05920\\
 & & \RG & 193.1 & 0.06747\\
\midrule
rough & $10^{-3}$ & Safe Landweber & 165.4 & 0.99581\\
 & & AA--Landweber & 35.4 & 0.91974\\
 & & \RG & 137.5 & 1.00106\\
\midrule
rough & $10^{-4}$ & Safe Landweber & 1445.5 & 0.57487\\
 & & AA--Landweber & 351.7 & 0.58201\\
 & & \RG & 999.0 & 0.57491\\
\bottomrule
\end{tabular}
\end{table}

AA--Landweber requires the fewest PDE solves in all four cases in \Cref{tab:1d-final}. \RG reduces work relative to safe Landweber with broadly comparable errors, but does not recover the full computational gain of unsafeguarded AA. For the rough truth at the lower noise level, its mean error is close to that of safe Landweber and smaller than that of AA--Landweber. Thus, in these tests, \RG retains part of the acceleration over safe Landweber while maintaining broadly comparable reconstruction errors.

\FloatBarrier
\subsection{Fixed-dimensional nonlinear rate benchmark}\label{sec:diag-rate}
To examine the benchmark rate in \Cref{thm:sqrtrate}, we consider a nonlinear diagonal model with an explicit source condition:
\begin{equation}\label{eq:diagmodel}
 \begin{gathered}
 F(x)_j=j^{-1}\left(x_j+\frac{\alpha}{2}x_j^2\right),\quad \alpha=0.1,\\
 F'(x)=Q_xT,\quad
 (Q_xv)_j=\frac{1+\alpha x_j}{1+\alpha x_j^\dagger}v_j,\quad T=F'(x^\dagger).
 \end{gathered}
\end{equation}
This factorization satisfies the range-invariance bounds in \Cref{ass:rate}.  The certified radius used in the localization audit satisfies $\kappa_0R=0.12<1$.
We choose
\[
 x_j^\dagger=0.55\sin(0.17j)j^{-0.55},\qquad
 \lambda_j^\dagger=
 \frac{1}{\sqrt{\zeta(1.01)}}(-1)^{j-1}j^{-0.505}.
\]
We set $x_0=x^\dagger-T^*\lambda^\dagger$. Both sequences belong to $\ell^2$, and the normalization gives $\|\lambda^\dagger\|_{\ell^2}=1$. All tests use the fixed truncation dimension $n=2048$.
For this setup, the exact-data localization check confirms $d_0+\Omega\sigma_0<R$, as required by \Cref{ass:margin}.

\begin{table}[H]
\centering
\renewcommand{\arraystretch}{1.00}
\caption{Mean iteration counts and reconstruction errors for the nonlinear diagonal problem. Slopes are fitted using all nine noise levels and the six lowest noise levels, respectively.}
\label{tab:diag-rate}
\footnotesize
\setlength{\tabcolsep}{3pt}
\begin{tabular}{lrrrrrr}
\toprule
& \multicolumn{2}{c}{Safe Landweber} & \multicolumn{2}{c}{AA--Landweber} & \multicolumn{2}{c}{\RG}\\
\cmidrule(lr){2-3}\cmidrule(lr){4-5}\cmidrule(lr){6-7}
$\delta_{\mathrm{rel}}$ & iterations & error & iterations & error & iterations & error\\
\midrule
$10^{-2}$ & 12.00 & $9.664\times10^{-3}$ & 8.00 & $5.613\times10^{-3}$ & 8.00 & $5.537\times10^{-3}$\\
$8\times10^{-3}$ & 12.45 & $8.857\times10^{-3}$ & 8.00 & $5.585\times10^{-3}$ & 8.00 & $5.510\times10^{-3}$\\
$5\times10^{-3}$ & 13.70 & $7.355\times10^{-3}$ & 9.00 & $5.387\times10^{-3}$ & 8.00 & $5.477\times10^{-3}$\\
$3\times10^{-3}$ & 18.30 & $5.637\times10^{-3}$ & 9.00 & $5.369\times10^{-3}$ & 8.00 & $5.461\times10^{-3}$\\
$10^{-3}$ & 34.85 & $3.264\times10^{-3}$ & 13.00 & $3.298\times10^{-3}$ & 13.00 & $3.338\times10^{-3}$\\
$8\times10^{-4}$ & 43.25 & $2.961\times10^{-3}$ & 16.00 & $2.939\times10^{-3}$ & 16.00 & $2.977\times10^{-3}$\\
$5\times10^{-4}$ & 55.75 & $2.328\times10^{-3}$ & 22.00 & $2.305\times10^{-3}$ & 25.05 & $2.423\times10^{-3}$\\
$3\times10^{-4}$ & 79.50 & $1.805\times10^{-3}$ & 30.00 & $1.839\times10^{-3}$ & 45.20 & $1.870\times10^{-3}$\\
$10^{-4}$ & 134.95 & $1.042\times10^{-3}$ & 58.80 & $1.082\times10^{-3}$ & 79.15 & $1.066\times10^{-3}$\\
\midrule
\multicolumn{1}{l}{Nine-level slope} & \multicolumn{2}{c}{0.487} & \multicolumn{2}{c}{0.363} & \multicolumn{2}{c}{0.360}\\
\multicolumn{1}{l}{Six-low-noise slope} & \multicolumn{2}{c}{0.496} & \multicolumn{2}{c}{0.473} & \multicolumn{2}{c}{0.480}\\
\bottomrule
\end{tabular}
\end{table}

Table~\ref{tab:diag-rate} shows that the fitted slopes at the six lowest noise levels are close to $1/2$. As shown in Figure~\ref{fig:diag-rate}, the two AA-based methods give smaller errors than safe Landweber at higher noise levels, but their errors vary only slightly across these levels. This plateau leads to smaller slopes when all nine noise levels are included in the fit. \RG achieves reconstruction errors comparable to AA--Landweber and requires fewer iterations than safe Landweber at every tested noise level. These results are consistent with the predicted low-noise convergence rate for \RG and show that it retains an iteration-count advantage over safe Landweber in this example.

\begin{figure}[H]
\centering
\includegraphics[width=.56\textwidth]{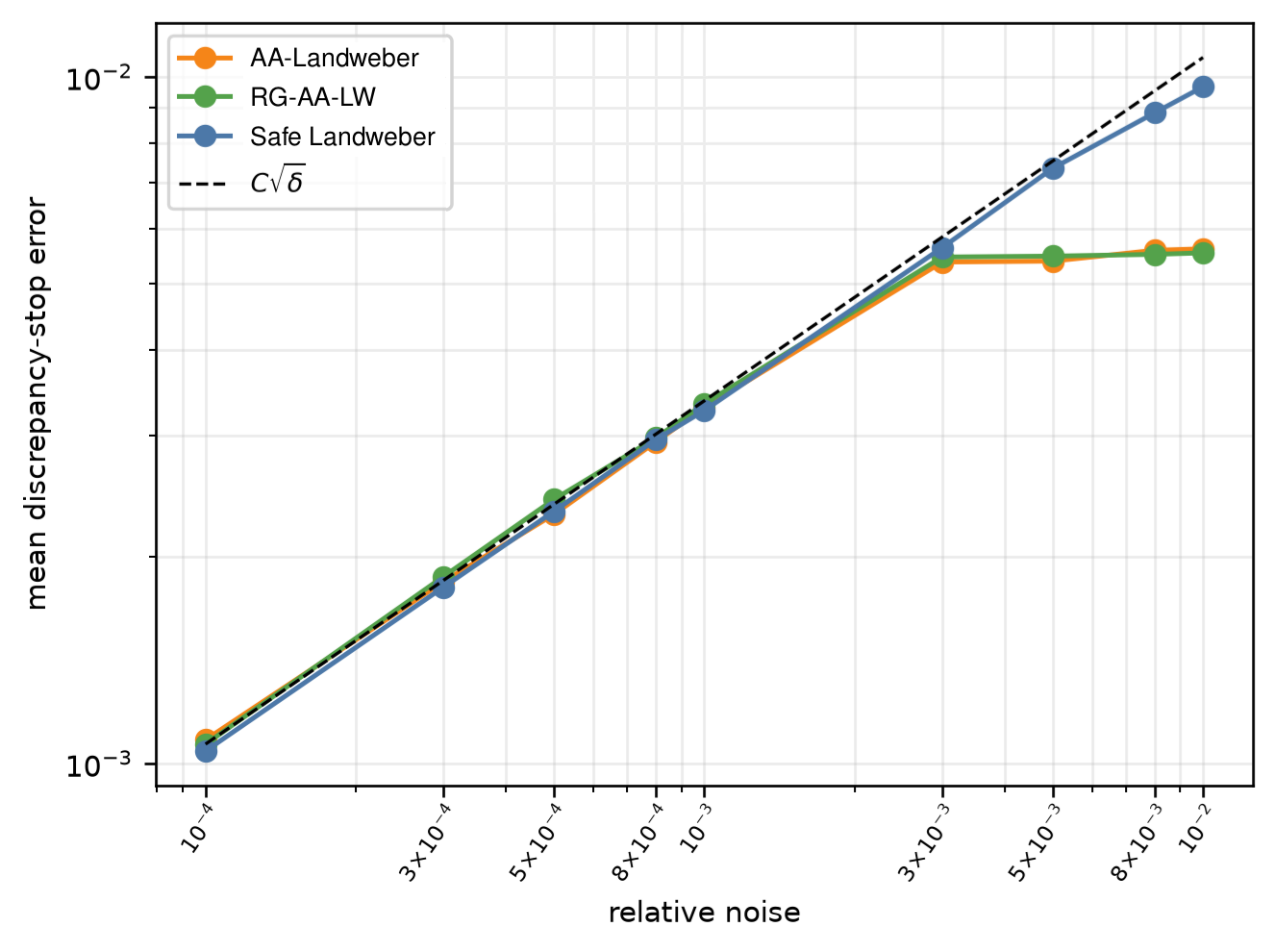}
\caption{Mean reconstruction errors for the nonlinear diagonal problem. The dashed line has slope $1/2$.}
\label{fig:diag-rate}
\end{figure}

We next compare the iteration histories of AA--Landweber and \RG for the same diagonal problem. The two methods use identical deterministic noisy data with $\delta_{\mathrm{rel}}=10^{-3}$ and the same proposal step $\mu=0.1t_0^\delta$, where $t_0^\delta$ is the initial certified step. Figure~\ref{fig:diag-aa-fail} shows one example of the residual norms and Anderson coefficient norms along the two trajectories.

In the left panel, AA--Landweber can sometimes exhibit pronounced residual oscillations, whereas \RG reduces the residual more steadily and reaches the discrepancy threshold much earlier. The right panel shows a larger early coefficient peak and more persistent coefficient fluctuations for AA--Landweber. \RG can also generate large coefficients, but these are not accompanied by comparable residual oscillations in this example. In \RG, the safeguard controls the correction actually accepted, rather than directly bounding the Anderson coefficients. The observed trajectories are consistent with this design and illustrate the more stable residual reduction achieved by \RG in this test.

\begin{figure}[H]
\centering
\includegraphics[width=.88\textwidth]{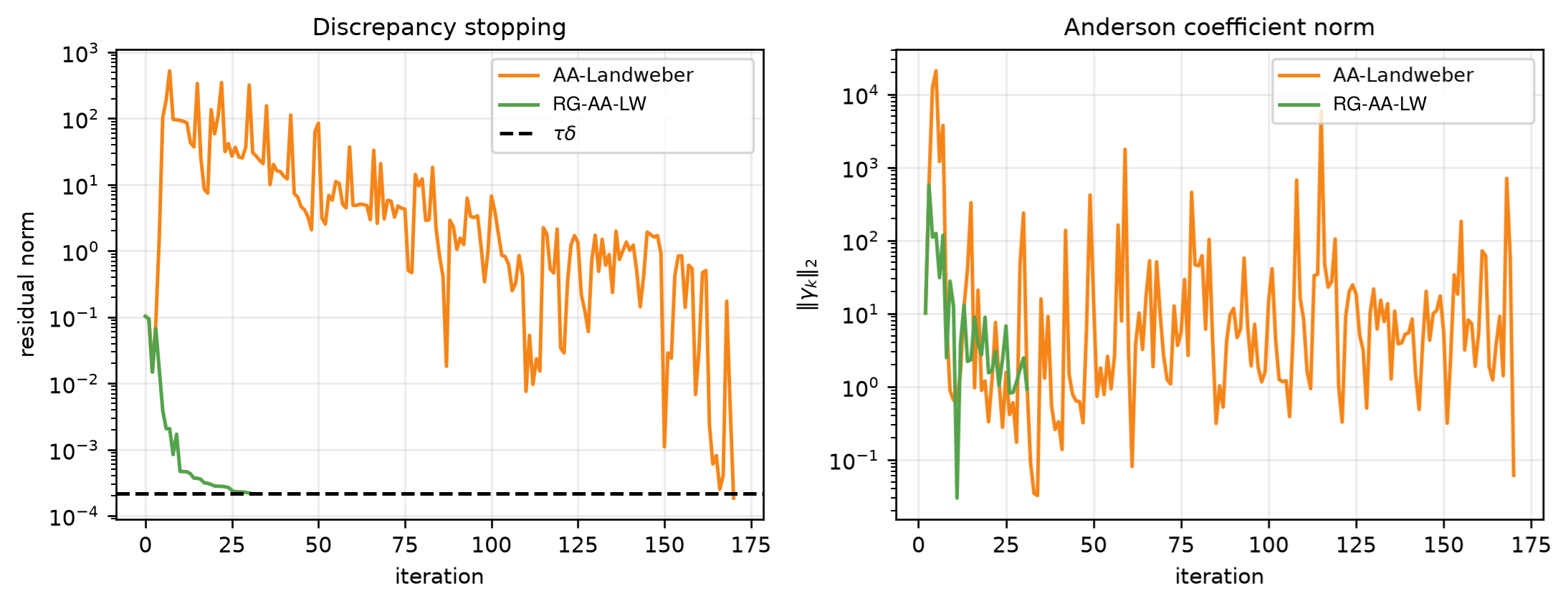}
\caption{Residual norms (left) and Anderson coefficient norms (right) for AA--Landweber and \RG. The dashed line marks the discrepancy threshold.}
\label{fig:diag-aa-fail}
\end{figure}

\FloatBarrier
\subsection{Two-dimensional potential inversion and strong baselines}
We next recover $q$ from internal observations of $u$ satisfying
\[
 -\Delta u+qu=1\quad\text{in }(0,1)^2,
 \qquad u=0\quad\text{on }\partial(0,1)^2.
\]
The PDE is discretized by standard five-point finite differences on a uniform grid with $256\times256$ interior points. We use $q_0\equiv0.8$ and $\delta_{\mathrm{rel}}=10^{-4}$ throughout this subsection. The smooth exact potential is
\begin{equation}\label{eq:2d-smooth-truth}
 \begin{aligned}
 q_{\mathrm{smooth}}^\dagger(x,y)
 ={}&0.8+\exp\!\left[-\frac{(x-0.35)^2+(y-0.35)^2}{0.025}\right]\\
 &+0.7\exp\!\left[-\frac{(x-0.70)^2+(y-0.65)^2}{0.035}\right].
 \end{aligned}
\end{equation}
The rough exact potential is
\begin{equation}\label{eq:2d-rough-truth}
 q_{\mathrm{rough}}^\dagger(x,y)=
 \begin{cases}
  1.8, & (x,y)\in R,\\
  2.4, & (x,y)\in D,\\
  1.25, & (x,y)\in B,\\
  0.7, & \text{otherwise},
 \end{cases}
\end{equation}
where
\[
 \begin{aligned}
 R&=(0.18,0.50)\times(0.20,0.68),\\
 D&=\{(x,y):(x-0.70)^2+(y-0.65)^2<0.12^2\},\\
 B&=(0.60,1)\times(0,0.35).
 \end{aligned}
\]
$R$, $D$, and $B$ are a rectangle, a disk, and a lower-right rectangle, respectively. Both truths are shown in \Cref{fig:2d-truths}.

To examine the effect of partial observation, we compare full-domain measurements with measurements on the first half of the grid rows in the $x$-direction and on the centered $128\times128$ block. These masked tests assess robustness beyond the verified range-invariance setting.

\begin{figure}[H]
\centering
\includegraphics[width=.90\textwidth]{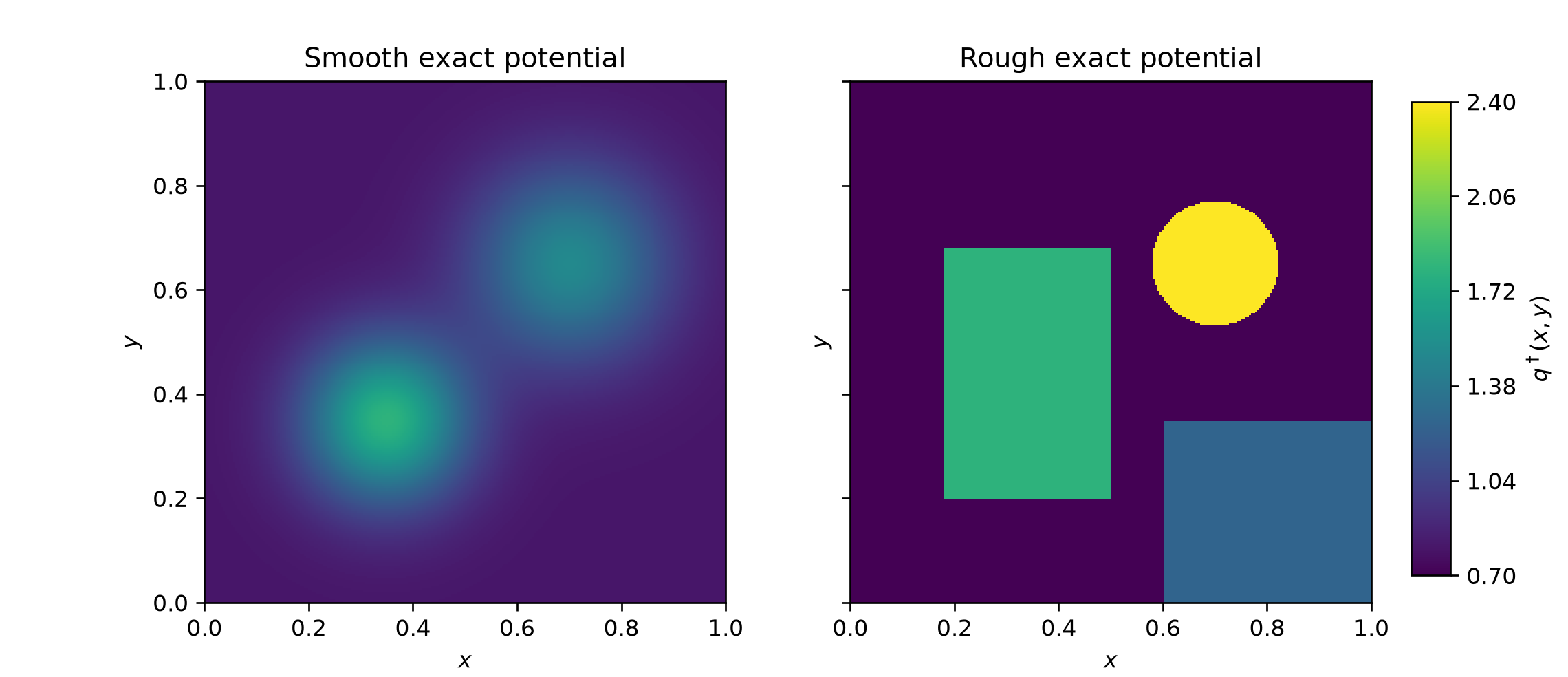}
\caption{$256^2$ exact potentials: smooth (left) and rough (right).}
\label{fig:2d-truths}
\end{figure}

\begin{table}[H]
\centering
\renewcommand{\arraystretch}{1.00}
\caption{Mean results for two-dimensional potential inversion under different observation regions at $\delta_{\mathrm{rel}}=10^{-4}$.}
\label{tab:2d-main}
\fontsize{9}{11}\selectfont
\setlength{\tabcolsep}{1.2pt}
\begin{tabular}{lrrrrrr}
\toprule
Case & \shortstack{AA--Landweber\\factorizations} & \shortstack{\RG\\factorizations} & \shortstack{AA--Landweber\\error} & \shortstack{\RG\\error} & \shortstack{AA--Landweber\\wall time (s)} & \shortstack{\RG{} wall\\time (s)}\\
\midrule
smooth, 100\% obs. & 16.00 & 14.00 & 0.02510 & 0.02320 & 10.27 & 8.66\\
smooth, 50\% obs. & 16.00 & 14.95 & 0.11891 & 0.11913 & 10.41 & 9.22\\
smooth, 25\% obs. & 9.00 & 10.75 & 0.04141 & 0.04123 & 5.49 & 6.25\\
rough, 100\% obs. & 118.95 & 82.90 & 0.21205 & 0.21063 & 85.70 & 56.40\\
rough, 50\% obs. & 37.15 & 32.70 & 0.38808 & 0.38795 & 27.77 & 23.31\\
rough, 25\% obs. & 57.85 & 60.40 & 0.28323 & 0.28266 & 40.87 & 40.15\\
\bottomrule
\end{tabular}
\end{table}

Table~\ref{tab:2d-main} shows that \RG requires fewer factorizations and less wall time than AA--Landweber under full and half-domain observation for both truths, with comparable reconstruction errors. This advantage does not extend to quarter-domain observation: slightly smaller errors require more factorizations, with longer wall time for the smooth truth and nearly unchanged time for the rough truth. Thus, reducing the observation region can change the computational benefit of the safeguard.

For comparison with stronger baselines, \Cref{tab:strong} includes Nesterov and IRGNM under half-domain observation. Nesterov uses the extrapolation parameter $\psi_k=(k-1)/(k+a-1)$ with $a=4$ \citep{HubmerRamlau2018}. For each truth, the IRGNM regularization schedule is selected on independent pilot data at $\delta_{\mathrm{rel}}=10^{-3}$, then fixed for testing at $\delta_{\mathrm{rel}}=10^{-4}$ with relative CG tolerance $10^{-5}$.

\begin{table}[H]
\centering
\renewcommand{\arraystretch}{1.00}
\caption{Mean results for the comparison with Nesterov and IRGNM under half-domain observation at $\delta_{\mathrm{rel}}=10^{-4}$.}
\label{tab:strong}
\footnotesize
\setlength{\tabcolsep}{3.5pt}
\begin{tabular}{llrrrrr}
\toprule
Truth & method & \shortstack{outer\\iterations} & factorizations & \shortstack{PDE\\solves} & error & \shortstack{wall time\\(s)}\\
\midrule
smooth & AA--Landweber & 15.00 & 16.00 & 31.0 & 0.11891 & 11.90\\
 & \RG & 13.95 & 14.95 & 28.9 & 0.11913 & 10.46\\
 & Nesterov & 34.05 & 35.05 & 69.1 & 0.11811 & 23.60\\
 & IRGNM & 12.00 & 13.00 & 223.0 & 0.11660 & 13.77\\
\midrule
rough & AA--Landweber & 36.15 & 37.15 & 73.3 & 0.38808 & 22.51\\
 & \RG & 31.70 & 32.70 & 64.4 & 0.38795 & 18.54\\
 & Nesterov & 71.90 & 72.90 & 144.8 & 0.38639 & 39.29\\
 & IRGNM & 14.00 & 15.00 & 347.3 & 0.38409 & 14.64\\
\bottomrule
\end{tabular}
\end{table}

All four methods satisfy the discrepancy criterion in all 20 realizations. \RG has the lowest PDE-solve count and wall time among the first-order methods, with errors close to those of AA--Landweber. IRGNM gives the smallest errors and uses the fewest outer iterations, but its inner iterations require substantially more PDE solves. It is faster than \RG for the rough truth and slower for the smooth truth. Nesterov slightly improves the errors of the AA-based methods at greater computational cost.

\FloatBarrier
\subsection{Mesh scaling for the two-dimensional potential problem}
We repeat the rough half-domain problem on grids from $32^2$ to $256^2$, keeping the remaining settings unchanged. This comparison examines whether the work reduction observed above persists as the number of unknowns increases.

\begin{table}[H]
\centering
\renewcommand{\arraystretch}{1.00}
\caption{Mean results under mesh refinement for the rough potential with half-domain observation.}
\label{tab:scaling}
\small
\begin{tabular}{llrrr}
\toprule
Grid & method & factorizations & wall time (s) & error\\
\midrule
$32^2$ & AA--Landweber & 49.00 & 0.227 & 0.38740\\
 & \RG & 29.65 & 0.125 & 0.38730\\
$64^2$ & AA--Landweber & 44.85 & 0.702 & 0.38851\\
 & \RG & 34.40 & 0.516 & 0.38839\\
$128^2$ & AA--Landweber & 44.25 & 6.597 & 0.38839\\
 & \RG & 30.05 & 4.424 & 0.38820\\
$256^2$ & AA--Landweber & 38.40 & 22.076 & 0.38808\\
 & \RG & 33.05 & 17.933 & 0.38803\\
\bottomrule
\end{tabular}
\end{table}

Table~\ref{tab:scaling} shows fewer factorizations and shorter wall times for \RG on every tested grid, while its reconstruction errors remain close to those of AA--Landweber. The computational benefit persists under mesh refinement, although the relative gain is smaller on the finest grid.

\FloatBarrier
\subsection{A semilinear PDE performance example}
Finally, we consider source inversion for a different nonlinear forward map:
\begin{equation}\label{eq:semilinear-source}
 -\Delta u+\kappa u^3=x\quad\text{in }(0,1)^2,
 \qquad u=0\quad\text{on }\partial(0,1)^2,
 \qquad F(x)=u(x).
\end{equation}
We use $\kappa=8$, with
\[
 F'(x)=(-\Delta+3\kappa u(x)^2)^{-1},
\]
and prescribe the six-mode exact state
\[
\begin{aligned}
u^\dagger(x,y)={}&0.204\sin(\pi x)\sin(\pi y)
 +0.090\sin(2\pi x)\sin(3\pi y)\\
&+0.066\sin(4\pi x)\sin(2\pi y)
 +0.048\sin(5\pi x)\sin(4\pi y)\\
&+0.033\sin(6\pi x)\sin(\pi y)
 +0.024\sin(7\pi x)\sin(6\pi y).
\end{aligned}
\]
The exact source and initial guess are
\[
 x^\dagger=-\Delta u^\dagger+\kappa(u^\dagger)^3,
 \qquad x_0=x^\dagger-F'(x^\dagger)^*\lambda^\dagger,
\]
where $\lambda^\dagger$ is a six-mode sine expansion normalized to $\|\lambda^\dagger\|=120$. Tests use a $256^2$ grid and $\delta_{\mathrm{rel}}=10^{-2}$, $10^{-3}$, and $10^{-4}$. The step length is selected on independent pilot data at $\delta_{\mathrm{rel}}=10^{-3}$ and then held fixed.

\begin{table}[H]
\centering
\renewcommand{\arraystretch}{1.00}
\caption{Semilinear source inversion: mean results.}
\label{tab:semilinear}
\small
\begin{tabular}{llrrr}
\toprule
$\delta_{\mathrm{rel}}$ & method & iterations & factorizations & error\\
\midrule
$10^{-2}$ & Safe Landweber & 8.00 & 31.00 & 0.23232\\
 & AA--Landweber & 4.00 & 19.00 & 0.21241\\
 & \RG & 4.00 & 19.00 & 0.21242\\
\midrule
$10^{-3}$ & Safe Landweber & 35.75 & 115.25 & 0.03159\\
 & AA--Landweber & 8.00 & 31.00 & 0.03527\\
 & \RG & 8.00 & 31.00 & 0.03528\\
\midrule
$10^{-4}$ & Safe Landweber & 73.45 & 228.35 & 0.00509\\
 & AA--Landweber & 24.00 & 69.00 & 0.00452\\
 & \RG & 25.00 & 72.00 & 0.00476\\
\bottomrule
\end{tabular}
\end{table}

\Cref{tab:semilinear} shows that the two AA-based methods require fewer iterations and factorizations than safe Landweber, while their reconstruction errors remain comparable. The recorded safeguard factors equal one throughout, so the outer safeguard never truncates an AA proposal. In this regime, the outer safeguard remains inactive, but the two methods still differ in their computation of the Anderson coefficient vector $\gamma_k$. The close agreement between the two AA-based methods therefore illustrates that \RG can recover AA behavior whenever the raw proposals already satisfy the correction limits; the remaining small differences are caused by the slightly different computation of $\gamma_k$.

\FloatBarrier
\endgroup
\section{Conclusion and future work}
We developed an Anderson-driven nonlinear Landweber method with explicit iterative-regularization guarantees.  Standard finite-memory AA generates the proposal, while a certified Landweber anchor and a summable, anchor-relative correction control what is actually committed.  Under the local tangential cone and localization assumptions, the exact-data trajectory converges strongly to a solution and the discrepancy-stopped noisy reconstructions converge to the same solution, without requiring local uniqueness.  Local uniqueness is needed only when this selected limit is to be identified with a prescribed solution.

Under the classical benchmark source condition, left range invariance, and a rate-compatible outer schedule, the fixed-reference ridge, finite memory, and residual-energy estimate provide source-element stability and the discrepancy-principle rate $O(\sqrt\delta)$.  The key point is that the proof accommodates the data-dependent Anderson coefficients directly through a controlled source-space perturbation.  The numerical results should be interpreted accordingly: the safeguard is not a uniform acceleration theorem.  It can preserve substantial AA gains in favorable and rough PDE regimes, can become nonintrusive when the raw proposal is already admissible, and can deliberately sacrifice some raw-AA speed when that is the price of maintaining the regularization mechanism.

The present analysis focuses on unprojected Hilbert-space iterations and the benchmark exponent $\nu=1/2$.  Natural next steps are rates for $0<\nu<1/2$ or variational source conditions, projection-active and nonsmooth constrained variants, and verifiable range-invariance surrogates for partial observations.  On the algorithmic side, a posteriori localization criteria, adaptive correction schedules, and an observable certification of the tangential-cone bound could reduce the conservatism of the sufficient margin $d_0+\Omega\sigma_0<R$ while retaining the regularization guarantee.  A broader abstraction of the present proposal--commitment principle to other iterative regularization schemes is another direction for future work.

\end{document}